\documentclass[a4paper]{article}
\usepackage{amssymb}
\newtheorem{Th}{Theorem}

\newtheorem{Lma}{Lemma}[section]
\newtheorem{Dfi}{Definition}

\newcommand{\be}{\begin{equation}}
\newcommand{\ee}{\end{equation}}
\newcommand{\R}{\mathbb{R}}
\newcommand{\N}{\mathbb{N}}
\newcommand{\C}{\mathbb{C}}

\newcommand{\reset}{\setcounter{equation}{0}\setcounter{Th}{0}\setcounter{Prop}{0}\setcounter{Co}{0}
\setcounter{Lm}{0}\setcounter{Rm}{0}}

\def\ti{\tilde}
\def\lf{\left}
\def\rg{\right}

\def\al{\alpha}
\def\la{\lambda}

\def\ep{\varepsilon}
\def\ds{\displaystyle}

\def\p{\partial}
\def\bn{\vec{n}}
\def\bh{\mathbf{h}}
\def\bbe{\vec{e}}
\def\bH{\vec{H}}

\def\bv{\vec{v}}
\def\bg{\vec{g}}
\def\bC{\vec{C}}
\def\bA{\vec{A}}
\def\bB{\vec{B}}
\def\bD{\vec{D}}
\def\bF{\vec{F}}
\def\bL{\vec{L}}
\def\bR{\vec{R}}
\def\bX{\vec{X}}
\def\bY{\vec{Y}}
\def\bw{\vec{w}}
\def\bbh{\vec{h}}
\def\bk{\vec{k}}
\def\br{\vec{r}}
\def\bc{\vec{c}}
\def\vp{\varphi}
\def\bp{\vec{\phi}}
\def\bep{\vec{\epsilon}}
\def\bc{\vec{c}}
\begin{document}
\title{Analysis aspects of Willmore surfaces}
\author{Tristan Rivi\`ere\footnote{Department of Mathematics, ETH Zentrum,
CH-8093 Z\"urich, Switzerland.}}
\maketitle
{\bf Abstract :} We found a new formulation to the Euler-Lagrange equation of the Willmore functional for immersed surfaces
in ${\R}^m$. This new formulation of Willmore equation appears to be of divergence form, moreover, the non-linearities
are made of jacobians. Additionally to that, if $\bH$ denotes the mean curvature vector of the surface, this new form 
writes ${\mathcal L}\bH=0$
where ${\mathcal L}$ is a well defined locally invertible self-adjoint operator.    These 3 facts  have numerous consequences 
in the analysis of Willmore surfaces.
One first consequence is that the long standing open problem to give a meaning to the Willmore Euler-Lagrange equation for immersions having only $L^2$ bounded second fundamental form  is now solved.   
We then establish the regularity of weak $W^{2,p}-$Willmore surfaces for any $p$ for which the Gauss map is continuous : $p>2$.
This is based on the proof of an $\epsilon-$regularity result for weak Willmore surfaces.
We establish then a weak compactness result for Willmore surfaces of energy less than $8\pi-\delta$ for every $\delta>0$.
This theorem is based on a point removability result we prove for Wilmore surfaces in ${\R}^m$. This result extends to arbitrary codimension the main result in \cite{KS3} established for surfaces in ${\R}^3$. Finally, we deduce from this point removability result the strong compactness,
modulo the M\"obius group action,
of Willmore tori below the energy level $8\pi-\delta$ in dimensions 3 and 4. The dimension 3 case was already solved in \cite{KS3}.
\section{Introduction}

Weak formulations of PDE offer not only the possibility to enlarge the class of solutions to the space of singular solutions but also provide a flexible setting
in which the analysis of smooth solutions becomes much more efficient. This is the idea that we want to illustrate in this paper by
introducing this new weak formulation of Willmore surfaces.

\medskip  

For a given oriented surface $\Sigma$ and a smooth positive immersion  $\Phi$ of $\Sigma$ into the Euclidian space ${\R}^m$, for some $m\ge 3$, we introduce first the Gauss map $\bn$ from $\Sigma$ into $Gr_{m-2}({\R}^m)$, the grassmanian of oriented $m-2-$planes of ${\R}^m$, which to every point $x$ in $\Sigma$ assigns  the unit $m-2$-unit vector defining the $m-2-$plane $N_{\Phi(x)}\Phi(\Sigma)$
orthogonal to the oriented tangent space $T_{\Phi(x)}\Phi(\Sigma)$. This map $\bn$ from $\Sigma$ into
$Gr_{m-2}({\R}^m)$ defines a projection map $\pi_{\bn}$  : for every vector $\xi$ in $T_{\Phi(x)}({\R}^m)$
$\pi_{\bn}(\xi)$ is the orthogonal projection of $\xi$ onto $N_{\Phi(x)}\Phi(\Sigma)$. Let then $\bB_x$ be the second fundamental form of the immersion $\Phi$ of $\Sigma$. $\bB_x$ is a symmetric bilinear form  on $T_x\Sigma$ with values into $N_{\Phi(x)}\Phi(\Sigma)$.
$\bB_x$ is given by $\bB_x=\pi_{\bn}\circ d^2\Phi$. By the mean of the ambiant scalar product in ${\R}^m$, which induces a metric
$g$ on $\Sigma$, we define the trace of $\bB_x$, $tr(\bB_x)$, which is a vector in $N_{\Phi(x)}\Phi(\Sigma)$
given by $tr(\bB_x)=\bB_x(e_1,e_1)+\bB_x(e_2,e_2)$ where $(e_1,e_2)$ is an arbitrary orthonormal basis of $T_x\Sigma$. The mean curvature vector $\bH(x)$ at $x$ of the immersion by $\Phi$ of $\Sigma$ is with theses notations the vector in $N_{\Phi(x)}\Phi(\Sigma)$  given by
\[
\bH(x)=\frac{1}{2}tr(\bB_x)\quad .
\]
In the case where $m=3$, $\bH(x)$ is the product of the mean value $H=1/2(\kappa_1+\kappa_2)$ of the principal curvatures $\kappa_1$, $\kappa_2$ of the surface at $\Phi(x)$ by $\bn$, the unit normal vector.

The so called {\it Willmore Functional} is then the following Lagrangian
\be
\label{0.1}
W(\Phi(\Sigma))=\int_{\Sigma}|\bH|^2\ d\,vol_g\quad.
\ee
where $d\,vol_g$ is the area form of the metric $g$ induced on $\Phi(\Sigma)$ by the canonical metric on ${\R}^m$.

\medskip

This lagrangian has been apparently first considered in the early 20th century in various works by Thomsen \cite{Tho}, Schadow
and a bit later by Blaschke \cite{Bla}. It has been reintroduced and more systematically studied in the framework of the conformal geometry of surfaces in space by Willmore in 1965 \cite{Wil}. Beyond conformal geometry this lagrangian plays an important role in various areas in science such as
molecular biology, where it has been considered as a surface energy for lipid bilayers
known as  {\it Helfrich Model} \cite{Hef}, such as non-linear elasticity in
solid mechanics  where it arises as limiting energy in thin plate theory (see \cite{FJM} for instance)
or even in general relativity where the lagrangian (\ref{0.1}) is the main term in the so called
{\it Hawking quasi local mass} (see \cite{Haw}, \cite{HI})...etc.
One of the reason for the genericity of this lagrangian is maybe the property discovered by White
\cite{Whi} for $m=3$ and proved by  B.Y. Chen \cite{Che} for arbitrary $m$ which says that the functional remains unchanged under the action of a conformal diffeomorphism of ${\R}^m$ (and even under conformal changes of metric of the ambiant space).

We are interrested in this work in the critical points of (\ref{0.1}) for perturbations of the form
$\Phi+t\phi$ where $\phi$ is an arbitrary smooth map from $\Sigma$ into ${\R}^m$. These critical points are the so called {\it Willmore surfaces}. Because of the invariance of the lagrangian
under the action of conformal transformations of ${\R}^m$, images of Willmore surfaces by such
conformal transformations are still Willmore. Examples of Willmore surfaces are minimal surfaces for which $\bH\equiv 0$ and which realize then absolut minimum of $W$. Willmore surfaces satisfy an Euler-Lagrange equation discovered by
Willmore for $m=3$ (though it was apparently known by it's predecessors on the subject Thomsen, Schadow and Blaschke in the twenties) and was established in it's full generality, for arbitrary $m$, by Weiner in \cite{Wei}.
Before presenting the equation we need the following notations :  Consider for every vector $\bw$ in $N_{\Phi(x)}\Phi(\Sigma)$ 
the symmetric endomorphism $A_x^{\bw}$ of $T_x\Sigma$ satisfying for every pair of vectors 
$\bX$ and  $\bY$ in $T_x\Sigma$ the identity $g(A_x^{\bw}(\bX),\bY)=B_x(\bX,\bY)\cdot \bw$,
where $\cdot$ denotes the standard scalar product in ${\R}^m$. The map which to $\bw$ assigns
the symmetric endomorphism $A^{\bw}_x$ of $T_x\Sigma$ for the scalar product $g$ is an homomorphism that we denote $A_x$ from $N_{\Phi(x)}\Phi(\Sigma)$ into $S_g\Sigma_x$,
the linear space of symmetric endomorphisms from $T_x\Sigma$ with respect to $g$. Denote
$\ti{A}_x$ the endomorphism of $N_{\phi(x)}\phi(\Sigma)$ obtained by composing the transpose
$^{t}A_x$ of $A_x$ with $A_x$ : $\ti{A}_x=$ $ ^{t}A_x\circ A_x$. Let $(\bbe_1,\bbe_2)$ be an orthonormal basis
of $T_x\Sigma$ and let $\bL$ be a vector in $N_{\Phi(x)}\Phi(\Sigma)$,
we have that $\ti{A}(\bL)=\sum_{i,j} \bB(\bbe_i,\bbe_j)\ \bB(\bbe_i,\bbe_j)\cdot\bL$

 With these notations, $\Phi$ is  a smooth 
Willmore immersion if and only if it solves the following Euler-Lagrange equation 
\be
\label{0.2} 
\Delta_\perp\bH-2|\bH|^2\ \bH+\ti{A}(\bH)=0\quad\quad,
\ee
where $\Delta_\perp$ is the negative covariant laplacian for the connection $D$ in the normal bundle $N\Phi(\Sigma)$ to $\Phi(\Sigma)$ issued from the ambiant scalar product in ${\R}^m$ :  for every 
section $\sigma$ of $N\Phi(\Sigma)$  one has $D_{\bX}\sigma:=\pi_{\bn}(\sigma_\ast\bX)$. In the particular case when $m=3$,  the mean curvature vector $\bH$
is oriented along the unit normal to $\Phi(\Sigma)$, $\bH=H\,\bn$,
and (\ref{0.2}) is equivalent to the following equation satisfied by the mean curvature function $H$ : 
\be
\label{0.3}
\Delta_g H+2H\ (|H|^2-K)=0\quad\quad,
\ee
where $\Delta_g$ is the negative laplace operator for the induced metric $g$ on $\Phi(\Sigma)$ and $K$ is the scalar curvature of $(\Sigma,g)$.

Despite their elegant aspects  equations (\ref{0.2}) and (\ref{0.3}) offer challenging mathematical difficulties. First of all one has to observe that the
highest order term  $\Delta_\perp\bH$ for (\ref{0.2}) or $\Delta_g\, H$ for (\ref{0.3}) is non-linear since  the metric $g$ defining 
the Laplace operator depends on the variable immersion $\Phi$. Another difficulty comes from the fact that  the Euler-Lagrange equations (\ref{0.2}), (\ref{0.3}) are a-priori non compatible with the Lagrangian
(\ref{0.1})  in the following sense.
Making the minimal regularity assumption which ensures that the Lagrangian (\ref{0.1}) is finite - 
the second fundamental form $\bB$ is $L^2$ on $\Phi(\Sigma)$ -
is not enough in order for the non-linearities in the equations (\ref{0.2}) or (\ref{0.3}) to have a distributional meaning : the expression $|\bH|^2\,\bH$
requires at least that $\bH$ is in $L^3$ and not only in $L^2$...etc.

Recently the author proved in \cite{Ri1} that  that any Euler Lagrange equation of any 2-dimensional conformally invariant lagrangian
with quadratic growth (such as the harmonic map equations into riemannian submanifolds or such as  the precribed mean curvature equation) can be written in divergence form.
 This divergence form has numerous consequences for the analysis of this equation. It permits 
 , in particular, to 
 extend the set of solutions to subspaces of distribution with very low regularity requirement.
 It seems that the analysis developped in \cite{Ri1} could be extended to other conformally invariant equations such as the harmonic map 
 equations into Lorentzian manifolds. Granting this observation together with the correspondance established by Bryant \cite{Bry} between
  Willmore surfaces in ${\R}^3$ and harmonic maps into the Minkowski sphere $S^{3,1}$ in ${\R}^{4,1}$, the author found not really the technic but at least a strong encouragement for looking for a divergence form to the Willmore Euler-Lagrange equation (\ref{0.1}). 
  
Our first main result in the present work is the following
 \begin{Th}
 \label{th-0.1}
 Willmore Euler-Lagrange Equation (\ref{0.2}) is equivalent to 
 \be
 \label{0.4}
 d\lf(\ast_g\, d\bH-3\ast_g\, \pi_{\bn}(d\bH)\rg)+d\lf(\bH\wedge d\bn\rg)=0
 \ee
 where $\ast_g$ is the Hodge operator on $\Sigma$ associated to the induced metric $g$ and where we make the implicit
 identification, by the mean of the standard volume form of ${\R}^m$, between $m-1$-vectors and vectors in ${\R}^m$
 (i.e. $\bH\wedge d\bn$ is then a 1 form on $\Sigma$ with values into ${\R}^m$).

 Assuming the immersion $\Phi$ is conformal from the flat disc $D^2=\Sigma$   into ${\R}^m$ then 
$\Phi$ is Willmore if and only if
\be
\label{0.5}
\Delta\bH-3\ div(\pi_{\bn}(\nabla \bH))-div\lf(\bH\wedge\nabla^\perp\bn\rg)=0
\ee
where the operators $\Delta$, $div$ and $\nabla$ are taken with respect to the flat metric on $D^2$
($\Delta=\p^2_{x_1}+\p^2_{x_2}$, $div X=tr\circ\nabla$ and $\nabla=(\p_{x_1},\p_{x_2})$). The operator $\nabla^\perp$ denotes the rotation by $\pi/2$ of $\nabla$ : $\nabla^\perp:=(-\p_{x_2},\p_{x_1})$.\hfill$\Box$
\end{Th}

This justifies then the following terminology : for a given map $\bn$ from $D^2$ into $G_{m-2}({\R}^m)$ we shall denote ${\mathcal L}_{\bn}$ and call it the {\it Willmore operator} the  operator which to a function $\bw$ from $D^2$ into ${\R}^m$ assigns
\be
\label{0.6}
{\mathcal L}_{\bn}\bw:=\Delta \bw-3\ div(\pi_{\bn}(\nabla\bw))-div(\bw\wedge\nabla\bn)\quad.
\ee
Although it was not difficult to  check it came as a nice surprise that this elliptic operator is self-adjoint : for any choice of 
map $\bn$ in $W^{1,2}(D^2,Gr_{m-2}({\R}^m))$ and for any choice of compactly supported maps $\bv$ and $\bw$ from $D^2$
into ${\R}^m$ we have
\be
\label{0.6a}
\int_{D^2}\bv\cdot{\mathcal L}_{\bn}\bw=\int_{D^2}{\mathcal L}_{\bn}\bv\cdot\bw\quad .
\ee
Another important information is the following 
 \[
\pi_{\bn}(\vec{v}):=(\vec{v}\cdot \bn(x))\cdot\bn(x)
\]
where $\cdot$ is the contraction operators between $p-$ and $q-$vectors producing $|p-q|$vectors in ${\R}^m$. 
Observe that 
\be
\label{0.7}
div(\pi_{\bn}(\nabla \bw))=\Delta\lf[(\bw\cdot \bn(x))\cdot\bn(x)\rg]-(\bw\cdot\nabla \bn(x))\cdot\bn(x)-(\bw\cdot \bn(x))\cdot\nabla\bn(x)
\ee
Thus, assuming now that the unit $m-2-$vector $\bn$ is in $W^{1,2}$,  the distribution
${\mathcal L}_{\bn}\bw$ given by (\ref{0.6}) is well defined for an arbitrary choice of $\bw$
in $L^2(D^2)$. This shows that the Euler-Lagrange equation in the form (\ref{0.4}) or (\ref{0.5})
is compatible with the Lagrangian (\ref{0.1}). Indeed the equation has a distributional sense
under the least possible regularity requirement for the immersion $\Phi(\Sigma)$ :
This minimal requirement is for the Gauss map to be in $W^{1,2}$ on $\Sigma$ with respect to the 
induced metric. 
\medskip

We assume first that $\Phi(\Sigma)$ is included in a graph $G_f$ of a map $f$ from $D^2$ into ${\R}^{m-2}$. In order to ensure that the Gauss map to that graph is $W^{1,2}$ on the graph, 
we are making the minimal regularity requirements on $f$ which are $f \in W^{2,2}$ for $m=3$ in
one hand and $f\in W^{1,\infty}\cap W^{2,2}$  for $m>3$ in the other hand.

We now define the notion of {\it Weak Willmore Graph}. To that purpose
it will be convenient to introduce locally a conformal parametrization of the graph $G_f$. 
The local existence of such a parametrization is given in \cite{MS} for $m=3$ and $f\in W^{2,2}(D^2)$
and by theorem 5.1.1 in \cite{Hel} for $m>3$ and   $f \in W^{1,\infty}\cap W^{2,2}$.

\begin{Dfi}
\label{df-0.2} [{\bf Weak Willmore Graphs}]
Let $f$ be a $W^{2,2}$ function from $D^2$ into ${\R}$ for $m=3$, respectively let $f$ be a Lipschitz and $W^{2,2}$ map
from $D^2$ into ${\R}^{m-2}$ for $m>3$, we say that the graph $G_f$ of $f$ is a Weak Willmore Graph when, in a conformal 
parametrization $\Phi $ from $D^2$ into $G_f$, the $L^2$ mean curvature vector $\bH$ and the $W^{1,2}$ Gauss map 
$\bn$ solve the equation .
\be
\label{0.8}
{\mathcal L}_{\bn}\bH=\Delta\bH-3\ div(\pi_{\bn}(\nabla \bH))-div\lf(\bH\wedge\nabla^\perp\bn\rg)=0\quad\quad\mbox{ in }{\mathcal D}'(D^2)
\quad.
\ee\hfill$\Box$
\end{Dfi}

Defininig $W^{2,2}$ Willmore immersion is a-priori problematic for the following reasons.
Let $\Phi$ be a $W^{2,p}$ map from a closed surface into ${\R}^m$. The assumption that $\Phi$
is an immersion means that $|\p_{x_1}\Phi\wedge\p_{x_2}\Phi|/|\nabla\Phi|^2>0$ a.e. Assuming that $p>2$
implies that $\nabla\Phi$ is continuous and therefore the immersion assumption 
gives a positive lower bound on $\Sigma$ to the function $|\p_{x_1}\Phi\wedge\p_{x_2}\Phi|/|\nabla\Phi|^2$. As a consequence we obtain that $\bn$ is also in $W^{1,p}$. Whereas this was a-priori not true for $p=2$ . The assumption that $p>2$ implies moreover 
that $\Phi(\Sigma)$ is made of a finite union of $W^{2,p}$ graphs because of the continuity of the Gauss map in that case. Hence,
 using definition~\ref{df-0.2}, the assumption
that the immersion $\Phi$ is in $W^{2,p}$ for some $p>2$ permits to define $W^{2,p}$ Willmore immersions.
Note that this extends to the situation where we only assume that the second derivatives of $\Phi$ are in the Lorentz space $L^{2,1}$.
This assumption on $\Phi$, which is weaker than the $W^{2,p}$ assumption for $p>2$, still ensures  the continuity of the Gauss map $\bn$ 
and that the second fundamental form of $\Phi(\Sigma)$ is in $L^2$ for the induced metric. It is a "border line assumption" in that respect.
$\nabla^2\Phi\in L^{2,1}(\Sigma)$  means that for some metric on $\Sigma$
\be
\label{0.9}
\int_0^{+\infty}\lf|\lf\{x\in\Sigma\ :\ |\nabla^2\Phi|(x)\ge\la\rg\}\rg|^\frac{1}{2}<+\infty\quad.
\ee
where $|\cdot|$ denotes the measure associated to the choice of metric we made on $\Sigma$.
An introduction to Lorentz spaces can be found in \cite{Ta2}.

\medskip

The second main result in the present work is the following
 
\begin{Th}
\label{th-0.3}
Let $\Sigma$ be a closed surface and let $\Phi$ be an immersion of $\Sigma$ in some euclidian space ${\R}^m$
 with 2 derivatives in the Lorentz space $L^{2,1}$.
Assume $\Phi$ is Willmore.
Then $\Phi(\Sigma)$ is the image of a real analytic immersion.\hfill$\Box$
\end{Th}

Note that since any $W^{2,p}$ map for $p>2$ has two derivatives in $L^{2,1}$, theorem~\ref{th-0.3} holds for $W^{2,p}$
Willmore immersions for $p>2$.

The previous theorem is a direct consequence of the following $\epsilon-$regularity result. 
This result was previously established
by Kuwert and Schaetzle in \cite{KS1} for strong solutions :  in the case where the surface is already assumed to be smooth.  

\begin{Th}
\label{th-III.1}{\bf [$\epsilon-$regularity for Weak Willmore Graphs.]}
Let $m$ be an integer larger than $2$. There exists $\ep>0$ such that the following holds. Let $f$ be a 
 function from the unit 2-dimensional
 disk $D^2$ into ${\R}^{m-2}$ with second derivatives in the Lorentz space $L^{2,1}$. Assume that the graph
 of $f$, $G_f\subset {\R}^m$, is Willmore (according to definition~\ref{df-0.2}).
Denote $\bn$ the Gauss map to this graph which is a map from $G_f$ into the grasmmanian of  unit $m-2-$vectors of ${\R}^m$.
 Then, under the following small energy assumption
\be
\label{III.1}
\int_{G_f}|\nabla \bn|^2\ d\,vol_{G_f}\le\ep\quad,
\ee
where the metric on $G_f$ is the one induced by the flat metric of ${\R}^m$, we have that the intersection of $G_f$ with
the cylinder $G_f({1/2}):=G_f\cap B^2_{1/2}(0)\times {\R}^{m-2}$ is $C^\infty$ and we have for every 
$k\in {\N}$
\be
\label{III.2}
\|\nabla^k \bn\|^2_{L^\infty(G_f(1/2))}\le C_k\ \int_{G_f}|\nabla \bn|^2\ d\,vol_{G_f}\quad,
\ee
where $C_k$ are  universal constants.   
\hfill$\Box$
\end{Th}
An interresting problem is to study whether this result extends or not  under the weaker assumption that $f$ is in $W^{2,2}$ for $m=3$
or respectively in $W^{2,2}\cap W^{1,\infty}$ for $m>3$. 

\medskip

Another fundamental question is to describe the "boundary" of the Moduli spaces 
of closed Willmore surfaces of given genus and bounded Willmore energy. In other words  one aims to describe the limiting behavior of sequences of Willmore surfaces $S_n$ with fixed topology and bounded Willmore energy. Modulo the action of the Moebius group of conformal 
transformations of 
${\R}^m$, which preseves Willmore Lagrangian, and therefore Willmore equation (\ref{0.4}), we can always fix the area of each $S_n$ to be equal to 1.  Now using Federer Fleming argument we can extract a subsequence to that sequence such that the current of integration on
$S_n$ converges for the Flat topology to some limiting integral current of integration  $S$
(see \cite{Fe} for the terminology of integral currents). Since $S_n$ has a uniformly bounded Willmore energy
and a fixed topology, the $L^2$ norm for the induced metric of it's second fundamental form  and hence the $W^{1,2}-$norm 
on the surface of the Gauss map are bounded. Applying then theorem~\ref{th-III.1} and a classical argument of concentration compactness
we then deduce that $S_n$ converges, in a suitable parametrization, in the $C^{k}$ topology to $S$ outside finitely many points $\{p_1,\cdots,p_k\}$. This strong convergence implies $S$ is a smooth Willmore surface
a-priori outside these points. The question to know whether these singular points are so called "removable" or not is then fundamental.
In the case where $S$  is a graph in a neighborhood of $p_k$  the regularity of $S$  about $p_k$  is given by  following result which extends to arbitrary codimensions the main result of \cite{KS3}.
\begin{Th}
\label{th-IV.1}
{\bf[Point removability for Willmore graphs.]}
Let $f$ be a continuous function from $B^2_1(0)$ into ${\R}^{m-2}$. Assume that $f$ realizes a $W^{2,p}-$Willmore graph over $D^2\setminus\{0\}$ for some $p>2$ and that the $W^{1,2}$ energy of the Gauss map on $G_f\setminus\{(0,f(0))\}$ is bounded.
Then the graph of $f$ over the whole disk $D^2$ is a $C^{1,\al}$ submanifold of ${\R}^m$ for every $\al<1$. Moreover, if $\bH$ denotes the mean curvature vector
of the graph, there exists a constant  vector $\bH_0$ such that $\bH(x)-\bH_0\log|x-x_0|$ is a $C^{0,\al}$ function on the graph, where $x_0=(0,f(0))$
and $|x-x_0|$ denotes the distance in the graph between $x$ and $x_0$.  If $\bH_0=0$ then $G_f$ is an analytic Willmore graph.
\hfill$\Box$
\end{Th}

Granting this result for $m=3$ Kuwert and Sch\"atzle were able to establish the fact that $S$ is a smooth Willmore surface in ${\R}^3$ under the assumption that the Willmore energy of the $S_n$ is less than $8\pi-\delta$ for any fixed $\delta>0$. This last fact ensures that
$S$ will be a graph about each $p_i$, $i=1\cdots k$ and that the residues $\bH_0$ will be equal to $0$ at each $p_i$. The arguments
 to prove
that, under the assumption that $W(S_n)<8\pi-\delta$, $S$ is a graph about each $p_i$ and that the residues $\bH_0$ at each $p_i$ are
equal to $0$ can be found page 344 of \cite{KS3} and are not specific to the codimension 1. Therefore combining them with our point removability result, theorem~\ref{th-IV.1}, with these arguments we can now state our last main result

\begin{Th}
\label{th-0.4} {\bf[Weak compactness of Willmore Surfaces below $8\pi$]}
Let $m$ be an arbitrary integer larger than $2$. Let $\delta>0$. Consider $S_n\subset{\R}^m$ to be a sequence of smooth closed Willmore  embeddings
with uniformly bounded topology, area equal to one and Willmore energy $W(S_n)$ bounded by $8\pi-\delta$. Assume that $S_n$
converges weakly as varifolds to some limit $S$ which realizes a non zero current. Then $S$ is a smooth Willmore embedding.
\hfill$\Box$
\end{Th}

Observe that the assumption that $S_n$ is a smooth Willmore immersion combined with the fact that $W(S_n)<8\pi$
implies that $S_n$ is an embedding due to a result by Li and Yau \cite{LY}.

Finally, combining again arguments  in \cite{KS3} (pages 350-351) together with our point removability result and a theorem
by Montiel in \cite{Mon} which in particular implies that any non-umbillic Willmore 2-spheres in ${\R}^4$ has Willmore
energy larger than $8\pi$ (this was known in ${\R}^3$ since the work of Bryant \cite{Bry}), we obtain the following.

\begin{Th}
\label{th-0.5} {\bf[Strong compactness of Willmore torii below $8\pi$]}
Let $m=3$ or $m=4$. Let $\delta>0$ arbitrary. The space of Willmore embedded torii in ${\R}^m$ having Willmore energy
less that $8\pi-\delta$ is compact up to M\"obius transformations under smooth convergence of compactly contained surfaces in 
${\R}^m$.\hfill$\Box$
\end{Th}

This extends to $m=4$ theorem 5.3 of \cite{KS3} where the above statement was proved for $m=3$.

\medskip

A crucial role in the analysis of Willmore surfaces that we develop in this paper is played by
the following observation.  We present it for $m=3$. Consider an $L^2$ map $\bw$ from $D^2$ into ${\R}^3$. Assume it is in the kernel of the Willmore operator 
${\mathcal L}_{\bn}$ for some map $\bn$ in $W^{1,2}(D^2,S^2)$ : ${\mathcal L}_{\bn}\bw=0$. 
Introduce the following Hodge decomposition
\be
\label{0.10}
\nabla\bw-3\pi_{N}(\nabla\bw)=\nabla \bA+\nabla^\perp \bB
\ee
for the boundary condition $\bA=0$ on $\p D^2$ for instance. Denote $(\bep_1,\bep_2,\bep_3)$ the canonical basis of ${\R}^3$. Then $\bA$ and $\bB=\sum_{i=1}^3B_i\bep_i$  solve the following 
equations
\be
\label{0.11}
\lf\{
\begin{array}{l}
\ds \Delta \bA=\nabla \bH\wedge\nabla^\perp \bn=\p_y\bH\wedge\p_x\bn-\p_x\bH\wedge\p_y\bn\\[5mm]
\ds \Delta B_i=3\,div(\pi_N(\nabla^\perp\bw))=3\sum_{j=1}^3\nabla(n_i\,n_j)\cdot \nabla^\perp w_j\\[5mm]
\ds\quad\quad=3\sum_{j=1}^3\p_y(n_i\,n_j)\ \p_x w_j-\p_x(n_i\,n_j)\ \p_y w_j
\end{array}
\rg.
\ee
where $\bw=\sum_{i=1}^3w_i\,\bep_i$ and $\bn=\sum_{i=1}^3n_i\,\bep_i$. The striking fact in the system (\ref{0.11}) is that all the non-linearities are linear combinations
of jacobians. This special algebraic structure plays a special role in geometric analysis. This was
probably first discovered by Wente in \cite{Wen}. An in-depth description of this is given in the book of H\'elein \cite{Hel}. This role is illustrated by the so called Wente type estimates (see section 3.1
of \cite{Hel}) which are used intensively in the present work.

\medskip
\noindent{\it Final Remarks : }
\begin{itemize}
\item[i)] The analysis we are developping in this work should give the direction for a new proof to Simon's result \cite{Si2} on the existence
of  embedded energy minimizing Willmore Tori in ${\R}^m$ for every $m\ge 3$.

\item[ii)] Our approach should be very useful in the study of the Willmore Flow initiated in the work of Kuwert and Sch\"atzle
 \cite{KS1}, \cite{KS2},\cite{KS3} (see also \cite{Sim}).

\item[iii)] Observe that for a Willmore surface the Hodge decomposition (\ref{0.10}) applied to the mean curvator
vector $\bH$ gives the following system.
\be
\label{0.12}
\lf\{
\begin{array}{l}
\ds \Delta \bA=\nabla \bH\wedge\nabla^\perp \bn\\[5mm]
\ds \Delta \bB=-3\nabla H\cdot\nabla^\perp\bn
\end{array}
\rg.
\ee
Since $|\nabla\bA|^2+|\nabla\bB|^2=4|\nabla H|^2+|H|^2\,|\nabla\bn|^2$, Wente estimates, with optimal Wente constants, applied to (\ref{0.12}) should give interresting lower bounds to
the Willmore energy of  Willmore immersions of Tori.

\item[iv)] Starting from the conservation law (\ref{0.4}) it is maybe possible to extend the notion 
of Willmore surfaces to varifolds with $L^2-$bounded generalized mean curvature (see \cite{Si1}
for this last notion). The need to extend Willmore surfaces to a larger class of object seems as natural
as the extension of the notion of minimal surfaces to minimal varifolds is.
\end{itemize}

\medskip

The paper is organised as follows. In section 2 we compute the conservation law satisfied by Willmore surfaces (proof of theorem~\ref{th-0.1}). In section 3 we give a proof of the $\ep-$regularity
for Willmore graphs, theorem~\ref{th-III.1}. In section 4 we give a proof of the point removability result
for Willmore graphs, theorem~\ref{th-IV.1}. A large part of the paper is devoted to the appendix in which we study,
outside the context of Willmore surfaces, various properties of the Willmore operator ${\mathcal L}_{\bn}$.

\medskip

\noindent{\it {\bf Acknowledgments} : The author is very grateful to Robert Bryant for stimulating discussions on Willmore surfaces
and in particular for having pointed out to him the work of Sebasti\`an Montiel quoted above.}

\section{Conservation laws for Willmore surfaces}
In this section the operators $\nabla$, $div$, $\Delta$ will be taken with respect
to the flat metric on the unit 2-dimensional disk $D^2=\{z\in{\C}\ ;\ |z|<1\}$
Let $Phi$ be a smooth conformal embedding of the unit disk $D^2$ into ${\R}^n$.
Denote $\Sigma=\Phi(D^2)$ and 
\[
e^{\la}=\lf|\frac{\p \Phi}{\p x_1}\rg|=\lf|\frac{\p \Phi}{\p x_2}\rg|
\]
Because of the
topology of $D^2$, the normal bundle to $\Sigma$ is trivial and there exists therefore
a smooth maps $\bn(z)=(\bn_1(z),\cdots,\bn_{m-2}(z))$ from $D^2$ into the orthonormal $m-2$ frames in 
${\R}^m$ satisfying
\[
(\bn_1(z),\cdots,\bn_{m-2}(z))\quad\mbox{ is a positive orthonormal basis to } N_{\Phi(z)}\Sigma\quad,
\]
where $N_{\Phi(z)}\Sigma=(T_{\Phi(z)}\Sigma)^\perp$ is the orthonormal $m-2$ plane to the tangent plane $T_{\Phi(z)}\Sigma$ of $\Sigma$ at $\Phi(z)$.
We denote by $(\bbe_1,\bbe_2)$ the orthonormal basis of $T_{\Phi(z)}\Sigma$ given by
\[
\bbe_i=e^{-\la}\ \frac{\p \Phi}{\p x_i}\quad .
\]
With these notations the second fundamental form $\bh$ which is a symmetric 2-form on $T_{\Phi(z)}\Sigma$  into  $N_{\Phi(z)}\Sigma$
is given by
\be
\label{II.1}
\bh=\sum_{\al,i,j}h^\al_{ij}\ \bn_\al\otimes(\bbe_i)^\ast\otimes(\bbe_j)^\ast\quad\mbox{ with }\quad h^\al_{ij}=-e^{-\la}\,\lf(\frac{\p \bn_\al}{\p x_i},\bbe_j\rg)
\ee
In particular the mean curvature vector $\bH$ is given by
\be
\label{II.2}
\bH=\sum_{\al=1}^{m-2} H^\al\,\bn_\al=\frac{1}{2}\sum_{\al=1}^{m-2}(h^\al_{11}+h^\al_{22})\, \bn_\al
\ee
Let $\bn$ be the $m-2$ vector of ${\R}^m$ given by $\bn=\bn_1\wedge\cdots\wedge\bn_2$. We identify vectors and $m-1$-vectors in ${\R}^m$
in the usual way. Hence we have for instance
\be
\label{II.0}
\bn\wedge \bbe_1=(-1)^{m-1} \bbe_2\quad\mbox{ and }\quad\bn\wedge \bbe_2=(-1)^{m-2} \bbe_1
\ee
Since $\bbe_1,\bbe_2,\bn_1\cdots\bn_{m-2}$ is a basis of $T_{\Phi(z)}{\R}^m$, we can write
for every $\al=1\cdots m-2$
\[
\nabla \bn_\al=\sum_{\beta=1}^{m-2}(\nabla \bn_\al,\bn_\beta)\, \bn_\beta+\sum_{i=1}^2(\nabla\bn_\al,\bbe_i)\,\bbe_i
\]
and consequently
\be
\label{II.3}
\bn\wedge\nabla^\perp\bn_\al=(-1)^{m-1}\ (\nabla^\perp \bn_\al,\bbe_1)\ \bbe_2
+(-1)^{m-2}\ (\nabla^\perp \bn_\al,\bbe_2)\ \bbe_1
\ee
Using the symmetry of the second fundamental form (i.e. $h^\al_{ij}=h^\al_{ji}$)
 and the conformality of $\Phi$ we have
\be
\label{II.4}
\lf(\frac{\p \bn_\al}{\p x_1},\bbe_2\rg)=\lf(\frac{\p \bn_\al}{\p x_2},\bbe_1\rg)\quad.
\ee
Thus, combining (\ref{II.3}) and (\ref{II.4}), we have
\be
\label{II.5}
\begin{array}{l}
\ds\bn\wedge\nabla^\perp\bn_\al\\[5mm]
\ds=(-1)^{m-1}\ \lf(
\begin{array}{c}
\ds-\lf(\frac{\p \bn_\al}{\p x_1},\bbe_2\rg)\\[3mm]
\ds\lf(\frac{\p \bn_\al}{\p x_1},\bbe_1\rg)
\end{array}
\rg)\ \bbe_2+(-1)^{m-1} \lf(
\begin{array}{c}
\ds\lf(\frac{\p \bn_\al}{\p x_2},\bbe_2\rg)\\[3mm]
\ds-\lf(\frac{\p \bn_\al}{\p x_2},\bbe_1\rg)
\end{array}
\rg)\ \bbe_1\\[5mm]
=(-1)^{m-1}\ \lf(
\begin{array}{c}
\ds-\lf(\frac{\p \bn_\al}{\p x_1},\bbe_2\rg)\\[3mm]
\ds-\lf(\frac{\p \bn_\al}{\p x_2},\bbe_2\rg)
\end{array}
\rg)\ \bbe_2+(-1)^{m-1} \lf(
\begin{array}{c}
\ds-\lf(\frac{\p \bn_\al}{\p x_1},\bbe_1\rg)\\[3mm]
\ds-\lf(\frac{\p \bn_\al}{\p x_2},\bbe_1\rg)
\end{array}
\rg)\ \bbe_1\\[5mm]
\ds +(-1)^{m-1}\lf[\lf(\frac{\p \bn_\al}{\p x_1},\bbe_1\rg)+\lf(\frac{\p \bn_\al}{\p x_2},\bbe_2\rg)\rg]\
\lf[\lf(
\begin{array}{c}
1\\[3mm]
0
\end{array}
\rg)\, \bbe_1+\lf(
\begin{array}{c}
0\\[3mm]
1
\end{array}
\rg)\, \bbe_2\rg]\quad .

\end{array}
\ee
This implies
\be
\label{II.6}
\begin{array}{l}
\ds\nabla \bn_\al+(-1)^{m-1} \bn\wedge\nabla^\perp\bn_\al\\[5mm]
\ds=\sum_{\beta=1}^{m-2}(\nabla \bn_\al,\bn_\beta)\ \bn_\beta-2e^{\la}\, H^\al\  \lf[\lf(
\begin{array}{c}
1\\[3mm]
0
\end{array}
\rg)\, \bbe_1+\lf(
\begin{array}{c}
0\\[3mm]
1
\end{array}
\rg)\, \bbe_2\rg]\quad .
\end{array}
\ee
Since $\p_x\Phi=e^{\la}\,\bbe_1$ and $\p_y\Phi=e^{\la}\, \bbe_2$ we finally get the following important identity
\be
\label{II.7}
\begin{array}{l}
\ds\nabla \bn_\al+(-1)^{m-1} \bn\wedge\nabla^\perp\bn_\al
\ds=\sum_{\beta=1}^{m-2}(\nabla \bn_\al,\bn_\beta)\ \bn_\beta-2H^\al\,\nabla\Phi\quad.
\end{array}
\ee  
Following the "Coulomb Gauge extraction Method" presented in the proof of lemma 4.1.3 of \cite{Hel}
we can choose a trivialization $(\bn_1\cdots\bn_{m-2})$ of the orthonormal frame bundle
associated to our trivial bundle $N\Sigma$
satisfying
\be
\label{II.8}
\forall 1\le\al,\beta\le m-2\quad\quad  div(\nabla \bn_\al,\bn_\beta)=0\quad.
\ee
Combining (\ref{II.7}) and (\ref{II.8}) we obtain in one hand
\be
\label{II.7a}
\begin{array}{l}
\ds div\lf(\nabla \bn_\al+(-1)^{m-1} \bn\wedge\nabla^\perp\bn_\al\rg)\\[5mm]
\ds=\sum_{\beta=1}^{m-2}(\nabla \bn_\al,\bn_\beta)\cdot\nabla \bn_\beta-2\nabla H^\al\cdot\nabla\Phi-2 H^\al\ \Delta\Phi\quad.
\end{array}
\ee  
In the other hand a classical elementary computation gives
\be
\label{II.8a}
\Delta\Phi=2\, e^{2\la}\ \bH\quad.
\ee
Therefore combining the two last identities we obtain
\be
\label{II.9}
\begin{array}{l}
\ds div\lf(\nabla \bn_\al+(-1)^{m-1} \bn\wedge\nabla^\perp\bn_\al\rg)\\[5mm]
\ds=\sum_{\beta=1}^{m-2}(\nabla \bn_\al,\bn_\beta)\cdot\nabla \bn_\beta-2\nabla H^\al\cdot\nabla\bp-4\, e^{2\la}\ H^\al\ \bH\quad.
\end{array}
\ee  
Multiplying this identity by $H^\al$, summing over $\al$ between $1$ and $m-2$ and using the expression of $H^\al\ \nabla \Phi$ given by (\ref{II.7}) we obtain
\be
\label{II.10}
\begin{array}{l}
\ds \sum_{\al=1}^{m-2}H^\al\ div\lf(\nabla \bn_\al+(-1)^{m-1} \bn\wedge\nabla^\perp\bn_\al\rg)
- \sum_{\al=1}^{m-2}\nabla H^\al\cdot\nabla \bn_\al\\[5mm]
\ds- (-1)^{m-1}\sum_{\al=1}^{m-2}
\nabla H^\al\cdot \bn\wedge\nabla^{\perp}\bn_\al\\[5mm]
\ds=\sum_{\al,\beta=1}^{m-2}H^\al\,(\nabla \bn_\al,\bn_\beta)\cdot\nabla \bn_\beta-\sum_{\al,\beta=1}^{m-2}\nabla H^\al\cdot(\nabla \bn_\al,\bn_\beta)\ \bn_\beta\\[7mm]
\ds \quad -4\, e^{2\la}\ H^\al\ \bH\quad.
\end{array}
\ee
For a purpose that will be clear later we rewrite equation (\ref{II.10}) in the following way
\be
\label{II.10a}
\begin{array}{l}
\ds \sum_{\al=1}^{m-2}H^\al\ \Delta\bn_\al- \sum_{\al=1}^{m-2}\nabla H^\al\cdot\nabla \bn_\al
-(-1)^{m-1}\sum_{\al=1}^{m-2}div \lf(H^\al\ \bn\wedge\nabla^\perp\bn_\al\rg)\\[5mm]
\ds=\sum_{\al,\beta=1}^{m-2}H^\al\,(\nabla \bn_\al,\bn_\beta)\cdot\nabla \bn_\beta-\sum_{\al,\beta=1}^{m-2}\nabla H^\al\cdot(\nabla \bn_\al,\bn_\beta)\ \bn_\beta\\[5mm]
\ds \quad-2(-1)^{m-1} \sum_{\al=1}^{m-2}H^\al\ div \lf(\bn\wedge\nabla^\perp\bn_\al\rg)-4\, e^{2\la}\ H^\al\ \bH\quad.
\end{array}
\ee
We shall give now an expression of $(-1)^{m-1} \sum_{\al=1}^{m-2}H^\al\ div \lf(\bn\wedge\nabla^\perp\bn_\al\rg)$ which will be usefull in the sequel. Using (\ref{II.3})
we have in one hand
\be
\label{II.11}
\begin{array}{l}
\ds (-1)^{m-1} div\lf( \bn\wedge\nabla^\perp\bn_\al\rg)\\[5mm]
\ds=(\nabla^\perp \bn_\al,\nabla\bbe_1)\ \bbe_2-(\nabla^\perp \bn_\al,\nabla\bbe_2)\ \bbe_1\\[5mm]
\ds+(\nabla^\perp \bn_\al,\bbe_1)\ \nabla\bbe_2-(\nabla^\perp \bn_\al,\bbe_2)\ \nabla\bbe_1\\[5mm]
\ds=\sum_{\beta=1}^{m-2}(\nabla^\perp\bn_\al,\bn_\beta)\ \lf[(\bn_\beta,\nabla\bbe_1)\ \bbe_2
-(\bn_\beta,\nabla \bbe_2)\ \bbe_1\rg]\\[5mm]
\ds+\sum_{\beta=1}^{m-2}\lf[(\nabla^\perp\bn_\al,\bbe_1)\ (\nabla\bbe_2,\bn_\beta)-(\nabla^\perp\bn_\al,\bbe_2)\ (\nabla\bbe_1,\bn_\beta)\rg]\ \bn_\beta\quad .
\end{array}
\ee
In the other hand, using the symmetry of $\bh$,
\be
\label{II.12}
\begin{array}{l}
\ds(\nabla^\perp\bn_\al,\bbe_1)\ (\nabla\bbe_2,\bn_\beta)-(\nabla^\perp\bn_\al,\bbe_2)\ (\nabla\bbe_1,\bn_\beta)\\[5mm]
\ds=\lf(\frac{\p \bbe_2}{\p x_1},\bn_\al\rg)\ \lf(\frac{\p \bbe_1}{\p x_2},\bn_\beta\rg)
-\lf(\frac{\p \bbe_2}{\p x_2},\bn_\al\rg)\ \lf(\frac{\p \bbe_1}{\p x_1},\bn_\beta\rg)\\[5mm]
\ds+\lf(\frac{\p \bbe_1}{\p x_2},\bn_\al\rg)\ \lf(\frac{\p \bbe_2}{\p x_1},\bn_\beta\rg)
-\lf(\frac{\p \bbe_1}{\p x_1},\bn_\al\rg)\ \lf(\frac{\p \bbe_2}{\p x_2},\bn_\beta\rg)\\[5mm]
\ds=e^{2\la}\lf[h^\al_{12}h^\beta_{12}-h^\al_{22}h^\beta_{11}+h^\al_{12}h^\beta_{12}-h^\al_{11}h^\beta_{22}\rg]\\[5mm]
\ds=e^{2\la}\ \sum_{i,j}h^\al_{ij}\,h^\beta_{ij}-4e^{2\la}\ H^\al\,H^\beta\quad.
\end{array}
\ee
Combining (\ref{II.11}) and (\ref{II.12}) we obtain
\be
\label{II.13}
\begin{array}{l}
\ds (-1)^{m-1} \sum_{\al=1}^{m-2}H^\al\ div\lf( \bn\wedge\nabla^\perp\bn_\al\rg)\\[5mm]
\ds=\sum_{\al,\beta=1}^{m-2}H^\al\ (\nabla^\perp\bn_\al,\bn_\beta)\ \lf[(\bn_\beta,\nabla\bbe_1)\ \bbe_2
-(\bn_\beta,\nabla \bbe_2)\ \bbe_1\rg]\\[5mm]
\ds\quad+e^{2\la}\ \sum_{\al,\beta,i,j}h^\al_{ij}\,h^\beta_{ij}-4e^{2\la}\ |\bH|^2\ \bH\quad .
\end{array}
\ee
That we rewrite in the following form
\be
\label{II.14}
\begin{array}{l}
\ds -2(-1)^{m-1} \sum_{\al=1}^{m-2}H^\al\ div\lf( \bn\wedge\nabla^\perp\bn_\al\rg)-4
e^{2\la}\ |\bH|^2\ \bH\\[5mm]
\ds=-2\sum_{\al,\beta=1}^{m-2}H^\al\ (\nabla^\perp\bn_\al,\bn_\beta)\ \lf[(\bn_\beta,\nabla\bbe_1)\ \bbe_2
-(\bn_\beta,\nabla \bbe_2)\ \bbe_1\rg]\\[5mm]
\ds\quad-2e^{2\la}\ \sum_{\al,\beta,i,j}h^\al_{ij}\,h^\beta_{ij}+4e^{2\la}\ |\bH|^2\ \bH\quad .
\end{array}
\ee
$\bH$ is a section of $N\Sigma$.  By definition the covariant (negative) laplacian of $\bH$
for the connection given by the orthogonal projection (w.r.t.  the standard scalar product in ${\R}^m$)
on the fibers is given by
\[
e^{2\la}\Delta_\perp\bH:=\pi_{\bn}div(\pi_{\bn}(\nabla \bH))\quad,
\]
where $\pi_{\bn}$ is the orthogonal projection on the fibers of $N\Sigma$. Using (\ref{II.8}) we have that
\be
\label{II.10b}
\begin{array}{l}
\ds e^{2\la}\Delta_\perp\bH:=\pi_{\bn}div\lf(\nabla H^\al\ \bn_\al+H^\al\ \lf(\nabla \bn_\al,\bn_\beta\rg)\,\bn_\beta\rg)\\[5mm]
\ds=\sum_\al\Delta H^\al\ \bn_\al
+2\sum_{\al,\beta}\nabla H^\al\cdot(\nabla\bn_\al,\bn_\beta)\,\bn_\beta\\[5mm]
\ds +\sum_{\al,\beta,\gamma} H^\al (\nabla\bn_\al,\bn_\beta)\ (\nabla\bn_\beta,\bn_\gamma)\ \bn_\gamma
\end{array}
\ee
Assuming now that our embedding is Willmore, it is equivalent to assume that $\bH$ satisfies (\ref{0.2})  (see \cite{Wei}) which means with the present notations that $\bh$ satisfies 
\be
\label{II.10c}
\Delta_\perp\bH+\sum_{i,j,\al,\beta}h^\al_{ij}\,h^\beta_{ij}\ H^\beta\ \bn_\al-2 |\bH|^2\ \bH=0\quad,
\ee
Thus, for a Willmore embedding the following identity holds
\be
\label{II.15}
\begin{array}{l}
\ds -2(-1)^{m-1} \sum_{\al=1}^{m-2}H^\al\ div\lf( \bn\wedge\nabla^\perp\bn_\al\rg)-4
e^{2\la}\ |\bH|^2\ \bH\\[5mm]
\ds=-2\sum_{\al,\beta=1}^{m-2}H^\al\ (\nabla^\perp\bn_\al,\bn_\beta)\ \lf[(\bn_\beta,\nabla\bbe_1)\ \bbe_2
-(\bn_\beta,\nabla \bbe_2)\ \bbe_1\rg]\\[5mm]
\ds\quad+2e^{2\la}\ \Delta_\perp\bH
\end{array}
\ee
Combining (\ref{II.10a}), (\ref{II.10b}) and (\ref{II.15}) we obtain that our conformal embedding $\Phi$
is Willmore if and only if the following identity holds
\be
\label{II.16}
\begin{array}{l}
\ds \sum_{\al=1}^{m-2}H^\al\ \Delta\bn_\al- \sum_{\al=1}^{m-2}\nabla H^\al\cdot\nabla \bn_\al
-(-1)^{m-1}\sum_{\al=1}^{m-2}div \lf(H^\al\ \bn\wedge\nabla^\perp\bn_\al\rg)\\[5mm]
\ds=\sum_{\al,\beta=1}^{m-2}H^\al\,(\nabla \bn_\al,\bn_\beta)\cdot\nabla \bn_\beta-\sum_{\al,\beta=1}^{m-2}\nabla H^\al\cdot(\nabla \bn_\al,\bn_\beta)\ \bn_\beta\\[5mm]
\ds-2\sum_{\al,\beta=1}^{m-2}H^\al\ (\nabla^\perp\bn_\al,\bn_\beta)\ \lf[(\bn_\beta,\nabla\bbe_1)\ \bbe_2
-(\bn_\beta,\nabla \bbe_2)\ \bbe_1\rg]\\[5mm]
\ds+2\sum_\al\Delta H^\al\ \bn_\al
+4\sum_{\al,\beta}\nabla H^\al\cdot(\nabla\bn_\al,\bn_\beta)\,\bn_\beta\\[5mm]
\ds +2\sum_{\al,\beta,\gamma} H^\al (\nabla\bn_\al,\bn_\beta)\ (\nabla\bn_\beta,\bn_\gamma)\ \bn_\gamma
\end{array}
\ee
We prove now that the previous identity (\ref{II.16}) can be written in divergence form. First we have
\be
\label{II.17}
\begin{array}{l}
\ds(\nabla^\perp\bn_\al,\bn_\beta)\ \lf[(\bn_\beta,\nabla\bbe_1)\ \bbe_2
-(\bn_\beta,\nabla \bbe_2)\ \bbe_1\rg]\\[5mm]
\ds=\lf[\lf(\frac{\p\bn_\al}{\p x_2},\bn_\beta\rg)\ \lf(\frac{\p\bn_\beta}{\p x_1},\bbe_1\rg)-\lf(\frac{\p\bn_\al}{\p x_1},\bn_\beta\rg)\ \lf(\frac{\p\bn_\beta}{\p x_2},\bbe_1\rg)\rg]\ \bbe_2\\[5mm]
\ds+\lf[\lf(\frac{\p\bn_\al}{\p x_2},\bn_\beta\rg)\ \lf(\frac{\p\bn_\beta}{\p x_1},\bbe_2\rg)-\lf(\frac{\p\bn_\al}{\p x_1},\bn_\beta\rg)\ \lf(\frac{\p\bn_\beta}{\p x_2},\bbe_2\rg)\rg]\ \bbe_1\\[5mm]
\ds=\lf[-2\lf(\frac{\p\bn_\al}{\p x_2},\bn_\beta\rg)\,e^{\la}\ H^\beta-\lf(\frac{\p\bn_\al}{\p x_2},\bn_\beta\rg)\ \lf(\frac{\p\bn_\beta}{\p x_2},\bbe_2\rg)\rg.\\[5mm]
\ds\lf.-\lf(\frac{\p\bn_\al}{\p x_1},\bn_\beta\rg)\ \lf(\frac{\p\bn_\beta}{\p x_1},\bbe_2\rg)\rg]\ \bbe_2\\[5mm]
\ds+\lf[-2\lf(\frac{\p\bn_\al}{\p x_1},\bn_\beta\rg)\,e^{\la}\ H^\beta-\lf(\frac{\p\bn_\al}{\p x_1},\bn_\beta\rg)\ \lf(\frac{\p\bn_\beta}{\p x_1},\bbe_1\rg)\rg.\\[5mm]
\ds\lf.-\lf(\frac{\p\bn_\al}{\p x_2},\bn_\beta\rg)\ \lf(\frac{\p\bn_\beta}{\p x_2},\bbe_1\rg)\rg]\ \bbe_1\\[5mm]
\ds=-2\ H^\beta\ (\nabla\bn_\al,\bn_\beta)\cdot\nabla\Phi\\[5mm]
\ds-(\nabla\bn_\al,\bn_\beta)\cdot\lf[
\lf(\nabla\bn_\beta,\bbe_2\rg)\ \bbe_2+\lf(\nabla\bn_\beta,\bbe_1\rg)\ \bbe_1\rg]
\end{array}
\ee
Thus we have 
\be
\label{II.18}
\begin{array}{l}
\ds-2\sum_{\al,\beta=1}^{m-2}H^\al\ (\nabla^\perp\bn_\al,\bn_\beta)\ \lf[(\bn_\beta,\nabla\bbe_1)\ \bbe_2
-(\bn_\beta,\nabla \bbe_2)\ \bbe_1\rg]\\[5mm]
\ds=4\sum_{\al=1}^{m-2}H^\al\  H^\beta\ (\nabla\bn_\al,\bn_\beta)\cdot\nabla\Phi\\[5mm]
\ds+2\sum_{\al,\beta=1}^{m-2}H^\al\ (\nabla\bn_\al,\bn_\beta)\cdot\lf[
\lf(\nabla\bn_\beta,\bbe_2\rg)\ \bbe_2+\lf(\nabla\bn_\beta,\bbe_1\rg)\ \bbe_1\rg]\quad .
\end{array}
\ee
Since $(\nabla\bn_\al,\bn_\beta)=-(\nabla\bn_\beta,\bn_\al)$, we deduce that
\be
\label{II.18a}
\sum_{\al,\beta=1}^{m-2}H^\al\  H^\beta\ (\nabla\bn_\al,\bn_\beta)\cdot\nabla\Phi\equiv 0\quad ,
\ee
and we finally get that
\be
\label{II.19}
\begin{array}{l}
\ds-2\sum_{\al,\beta=1}^{m-2}H^\al\ (\nabla^\perp\bn_\al,\bn_\beta)\cdot\lf[(\bn_\beta,\nabla\bbe_1)\ \bbe_2
-(\bn_\beta,\nabla \bbe_2)\ \bbe_1\rg]\\[5mm]
\ds+2\sum_{\al,\beta,\gamma}H^\al\ (\nabla\bn_\al,\bn_\beta)\cdot(\nabla\bn_\beta,\bn_\gamma)\ \bn_\gamma\\[5mm]
\ds=2\sum_{\al,\beta=1}^{m-2}H^\al\ (\nabla\bn_\al,\bn_\beta)\cdot\nabla\bn_\beta\quad .
\end{array}
\ee
Substituting this identity in (\ref{II.16}) we obtain that the conformal embedding $\Phi$ is Willmore
if and only if
\be
\label{II.20}
\begin{array}{l}
\ds \sum_{\al=1}^{m-2}H^\al\ \Delta\bn_\al- \sum_{\al=1}^{m-2}\nabla H^\al\cdot\nabla \bn_\al
-(-1)^{m-1}\sum_{\al=1}^{m-2}div \lf(H^\al\ \bn\wedge\nabla^\perp\bn_\al\rg)\\[5mm]
\ds=3\sum_{\al,\beta=1}^{m-2}H^\al\,(\nabla \bn_\al,\bn_\beta)\cdot\nabla \bn_\beta+2\sum_{\al=1}^{m-2}\Delta H^\al\ \bn_\al\\[5mm]
\ds+3\sum_{\al,\beta=1}^{m-2}\nabla H^\al\cdot(\nabla\bn_\al,\bn_\beta)\,\bn_\beta
\end{array}
\ee
Using the condition (\ref{II.8}) satisfied by our special choice of trivialization of the normal bundle,
 we rewrite the previous identity in the following way 
\be
\label{II.21}
\begin{array}{l}
\ds \sum_{\al=1}^{m-2}H^\al\ \Delta\bn_\al- \sum_{\al=1}^{m-2}\nabla H^\al\cdot\nabla \bn_\al
-(-1)^{m-1}\sum_{\al=1}^{m-2}div \lf(H^\al\ \bn\wedge\nabla^\perp\bn_\al\rg)\\[5mm]
\ds=3\ div\lf[\sum_{\al,\beta=1}^{m-2}H^\al\,(\nabla \bn_\al,\bn_\beta)\ \bn_\beta\rg]+2\sum_{\al=1}^{m-2}\Delta H^\al\ \bn_\al\quad.
\end{array}
\ee
Observe in one hand that 
\be
\label{II.22}
-(-1)^{m-1}\sum_{\al=1}^{m-2}H^\al\ \bn\wedge\nabla^\perp\bn_\al=
\sum_{\al=1}^{m-2}H^\al\ \nabla^\perp\bn_\al\wedge\bn=\nabla^\perp\bH\wedge\bn\quad,
\ee
and in the other hand that
\be
\label{II.23}
\sum_{\al,\beta=1}^{m-2}H^\al\,(\nabla \bn_\al,\bn_\beta)\ \bn_\beta=\sum_{\beta=1}^{m-2}\,
(\nabla\bH,\bn_\beta)\ \bn_\beta\quad.
\ee
Using these 2 observations the Willmore equation (\ref{II.21}) becomes
\be
\label{II.24}
\begin{array}{l}
\ds\Delta\bH-3\sum_{\al=1}^{m-2}\nabla H^\al\cdot\nabla\bn_\al-3\sum_{\al=1}^{m-2}\Delta H^\al\ \bn_\al
+div\lf(\nabla^\perp\bH\wedge\bn\rg)\\[5mm]
\ds=3\,div\lf(\sum_{\beta=1}^{m-2}(\nabla\bH,\bn_\beta)\ \bn_\beta\rg)\quad.
\end{array}
\ee
Observe now that 
\be
\label{II.25}
\begin{array}{l}
\ds\sum_{\al=1}^{m-2}\nabla H^\al\ \bn_\al+\sum_{\beta=1}^{m-2}(\nabla\bH,\bn_\beta)\ \bn_\beta\\[5mm]
\ds=\nabla\bH-\sum_{\al=1}^{m-2} H^\al\ \nabla\bn_\al+\sum_{\al,\beta=1}^{m-2} H^\al\ (\nabla\bn_\al,\bn_\beta)\ \bn_\beta\\[5mm]
\ds= \nabla\bH-\sum_{\al=1}^{m-2} H^\al\ (\nabla\bn_\al,\bbe_1)\ \bbe_1-\sum_{\al=1}^{m-2} H^\al\ (\nabla\bn_\al,\bbe_2)\ \bbe_2\\[5mm]
\ds=\nabla\bH-(\nabla\bH,\bbe_1)\ \bbe_1-(\nabla\bH,\bbe_2)\ \bbe_2=\pi_{\bn}(\nabla\bH)
\end{array}
\ee
So finally we have proved that the Willmore equation can be written in the following form
\be
\label{II.26}
\Delta\bH-3\ div(\pi_{\bn}(\nabla \bH))+div\lf(\nabla^\perp\bH\wedge\bn\rg)=0
\ee
\section{$\epsilon-$regularity for Weak Willmore graphs}

This section is devoted to the proof of the $\ep-$regularity theorem~\ref{th-III.1}

\medskip

Let then $f$ be a function from $D^2$ into ${\R}^{m-2}$ with 2 derivatives in the Lorentz space $L^{2,1}$ realizing a Willmore graph.
We consider the bilipschitz conformal parametrization from $D^2$ into $G_f$ given by theorem 5.1.1 of \cite{Hel}. 
Call this parametrization $\Phi$.
According to definition~\ref{df-0.2} this means that the mean curvator $\bH$ satisfies the equation (\ref{0.8}).
Denote by $\chi$  a smooth cut-off function equal to 1 on $D^2_{1/2}$ and compactly supported in $D^2$. Since ${\mathcal{L}}_{\bn} \bH=0,$ we have
\be
\label{III.c46}
\begin{array}{rl}
\ds{\mathcal{L}}_{\bn}(\chi\, \bH)=&\ds2 div(\nabla\chi\ \bH)-\bH\,\Delta\chi-6\, div(\pi_{\bn}(\bH)\,\nabla\chi)\\[5mm]
 &\ds+3(\bH\cdot\nabla\bn)\cdot\bn.\nabla\chi+3(\bH\cdot\bn)\cdot\nabla\bn.\nabla\chi\\[5mm]
&\ds-3\Delta\chi\pi_{\bn}(\bH)-\bH\wedge\nabla^\perp\bn.\nabla\chi
\end{array}
\ee
Let 
\be
\label{III.c47}
\bg_1=2 div(\nabla\chi\ \bH)-\bH\,\Delta\chi-6\, div(\pi_{\bn}(\bH)\,\nabla\chi)-3\Delta\chi\pi_{\bn}(\bH)
\ee
and
\be
\label{III.c48}
\bg_2=3(\bH\cdot\nabla\bn)\cdot\bn.\nabla\chi+3(\bH\cdot\bn)\cdot\nabla\bn.\nabla\chi
-\bH\wedge\nabla^\perp\bn.\nabla\chi
\ee
We have
\be
\label{III.c49}
\|\bg_1\|^2_{H^{-1}(D^2)}\le C\int_{B_1\setminus B_{1/2}}|\bH|^2
\ee
and
\be
\label{III.c50}
\|\bg_2\|_{L^1(D^2)}\le C\int_{D^2\setminus D^2_{1/2}}|\bH|\,|\nabla \bn|\quad.
\ee
Denote $\bv_1$ the solution of (\ref{III.7}) given by lemma~\ref{lm-III.1} for $\bg=\bg_1$ and let
$\bv_2$ be the solution of (\ref{III.32}) given by lemma~\ref{lm-III.3} for $\bg=\bg_2$.
We have in particular
\be
\label{III.c51}
\|\nabla\bv_1\|_{L^2}\le C\lf[\int_{D^2\setminus D^2_{1/2}}|\bH|^2\rg]^\frac{1}{2}\quad,
\ee
and
\be
\label{III.c52}
\|\nabla\bv_2\|_{L^{2,\infty}}\le C\int_{D^2\setminus D^2_{1/2}}|\bH|\,|\nabla \bn|\le
C\ep^\frac{1}{2}\lf[\int_{D^2\setminus D^2_{1/2}}|\bH|^2\rg]^\frac{1}{2}\quad.
\ee
$\bv:=\chi \bH-\bv_1-\bv_2$ is in $L^2(D^2)$ and satisfies ${\mathcal{L}}_{\bn}\bv=0$
we also observe that since $\chi$ is compactly supported in $D^2$, $\nabla\bv$ is
a sum of a compactly supported distribution in the interior of $D^2$ and a $L^{2,\infty}$
function. The trace of $\bv$ on $\p D^2$ is therefore well defined and equals  zero.  Assuming $\nabla \bn$
is in the Lorentz space $L^{2,1}(D^2)$  , we are now 
in position to apply lemma~\ref{lma-A.8} and we deduce that $\bv$ is identically $0$.
Thus we have proved that $\nabla(\chi\bH)$ is $L^{2,\infty}$ and
\be
\label{III.c53}
\|\nabla (\chi\bH)\|_{L^{2,\infty}(D^2)}\le C\lf[\int_{D^2\setminus D^2_{1/2}}|\bH|^2\rg]^\frac{1}{2}\quad.
\ee
We have ${\mathcal L}_{\bn}(\chi\bH)={\mathcal L}_{\bn}(\bH)=0$ on $D^2_{1/2}$. We proceed to the following
 Hodge decomposition of
 $\nabla \bH-3\pi_{\bn}(\nabla\bH)$
on $D^2_{1/2}$ : $\nabla \bH-3\pi_{\bn}(\nabla\bH)=\nabla C+\nabla^\perp D+\br$ with the boundary
conditions $C=0$ on $\p D^2_{1/2}$ and $\p D/\p\nu=0$ on $D^2_{1/2}$ and where $\br$ is harmonic.
 The following equations then holds
In one hand
\be
\label{III.c54}
\lf\{
\begin{array}{l}
\ds\Delta C=div(\bH\wedge\nabla^\perp\bn)\quad\quad\mbox{ in }D^2_{1/2}\\[5mm]
\ds C=0\quad\quad\mbox{ on }\p D^2_{1/2}
\end{array}
\rg.
\ee
In the other hand
\be
\label{III.c55}
\lf\{
\begin{array}{l}
\ds\Delta D=3\, div(\pi_{\bn}(\nabla^\perp \bH))\quad\quad\mbox{ in }D^2_{1/2}\\[5mm]
\ds\frac{\p D}{\p \nu}=0\quad\quad\mbox{ on }\p D^2_{1/2}
\end{array}
\rg.
\ee
Using now that the right-hand-sides of (\ref{III.c54}), (\ref{III.c55}) are jacobians of $\bH$ and $\bn$,
since $\nabla\bn\in L^2(D^2_{1/2})$ and $\nabla\bH\in L^{2,\infty}(D^2_{1/2})$, we have, using the Wente estimate (\ref{III.38}),
\be
\label{III.c56}
\|\nabla C\|_{L^2(D^2_{1/2})}+\|\nabla D\|_{L^2(D^2_{1/2})}\le C\|\nabla\bn\|_{L^2(D^2_{1/2})}\ \|\nabla\bH\|_{L^{2,\infty}(D^2_{1/2})}
\ee
Since $\br$ is harmonic we have that $\|\br\|_{L^2(D_{1/4})}\le C\|\br\|_{L^{2,\infty}(D_{1/2})}$. Combining this last inequality together with
(\ref{III.c56}), and the fact that $| \nabla \bH-3\pi_{\bn}(\nabla\bH)|\ge|\nabla \bH|$, we have established in particular that
\be
\label{III.c57}
\|\nabla \bH\|_{L^2(D^2_{1/4})}\le C\|\nabla\bn\|_{L^2(D^2_{1/2})}\ \|\nabla\bH\|_{L^{2,\infty}(D^2_{1/2})}+\|\nabla\bH\|_{L^{2,\infty}(D^2_{1/2})}
\ee
Using one more time the same Hodge decomposition but on $D^2_{1/4}$ instead of $D^2_{1/2}$, and replacing Wente's estimate
(\ref{III.38}) by Wente's inequality  (3.47) in theorem 3.4.1 of \cite{Hel} (which was originally obtained
by Luc Tartar in \cite{Ta1}) and argueing similarly as before we get
\be
\label{III.c58}
\|\nabla \bH\|_{L^{2,1}(D^2_{1/8})}\le C\|\nabla\bn\|_{L^2(D^2_{1/4})}\ \|\nabla\bH\|_{L^{2}(D^2_{1/4})}+\|\nabla\bH\|_{L^{2,\infty}(D^2_{1/4})}
\ee
Using now equation (\ref{II.8}), since $\nabla \la\in L^{2,1}(D^2)$ (see \cite{Hel}), we have that $\nabla(\Delta\bp)$ is in
$L^{2,1}(D^2_{1/8})$ and this implies that $\nabla\bn$ is in $L^\infty(D^2_{1/8})$.  Moreover, combining (\ref{III.c53}), (\ref{III.c57})
and (\ref{III.c58}) we obtain
\be
\label{III.c59}
\|\nabla\bn\|_{L^\infty(D^2_{1/8})}^2\le C\ \int_{D^2}|\nabla\bn|^2\quad,
\ee
and theorem~\ref{th-III.1} is proved.\hfill$\Box$
 
\section{Point removability for Willmore graphs.}

This section is devoted to the proof of theorem~\ref{th-IV.1}

Under the assumptions of the theorem we shall consider  the Lipshitz conformal parametrization $\Phi$ from $D^2$ into $G_f$, the graph of $f$, obtained by
following the arguments of \cite{KS3} (pages 332-334) which are based on the use of Huber's result on conformal parametrizations of complete surfaces in ${\R}^m$ 
\cite{Hub}
together with the estimates given by the work of   M\"uller and Sver\'ak \cite{MS}     . The preimage of $(0,f(0))$ by $\Phi$ being
in the inside of $D^2$ we can always assume (after a possible composition by a Moebius transformation of $D^2$) that $\Phi^{-1}((0,f(0))=0$.
Using theorem~\ref{th-III.1} we have that $\Phi$ is $C^\infty$ in $D^2\setminus\{0\}$, moreover, from (\ref{III.2}), there exist a positive function $\delta(r)$ 
going to zero as $r$ goes to zero such that
\be
\label{IV.6}
\forall x\in D^2\setminus\{0\}\quad\quad |x|\,|\nabla \bn(x)|+|x|^2\,|\nabla^2\bn|\le\delta(|x|)\quad.
\ee
The distances are taken with respect to the flat metric on $D^2$, it is however equivalent to the distance 
for the induced metric in the graph due to the estimates in \cite{MS}.

A positive number $\ep$ being given, we can restrict to a smaller disk and dilate to ensure that
\be
\label{IV.7}
\||x|\,|\nabla \bn|(x)\|^2_\infty+\int_{D^2}|\nabla\bn|^2\le\ep\quad.
\ee
Since $\bH$ is in $L^2(D^2)$ the distribution ${\mathcal L}_{\bn}\bH$ makes sense in ${\mathcal D}'(D^2)$. Moreover, since $\Phi$ is Willmore in $D^2\setminus\{0\}$, the distribution ${\mathcal L}_{\bn}\bH$ is supported in zero and is therefore a finite linear combination of derivatives of the Dirac
mass at the origin. Since ${\mathcal L}_{\bn}\bH$ is a sum of an $H^{-2}$ distribution and derivatives
of $L^1$ functions it can only be proportional to the Dirac mass centered at the origin itself :
\be
\label{IV.8}
{\mathcal L}_{\bn}\bH= \bc_0\ \delta_0\quad .
\ee
In anticipation to the result we introduce the constant $\bH_0$ satisfying $\bc_0= -4\pi\bH_0$.
Let $\bL$ the solution to the following problem given by lemma~\ref{lm-III.4}
\be
\label{IV.9}
\lf\{
\begin{array}{l}
\ds{\mathcal L}_{\bn}\bL=-4\pi\bH_0\quad\quad\mbox{ in }D^2\quad ,\\[5mm]
\ds \bL=0\quad\quad\mbox{ on }\p D^2\quad.
\end{array}
\rg.
\ee
We have $\nabla\bL\in L^{2,\infty}$. Since $\bn$ is smooth in $D^2\setminus\{0\}$ and satisfy 
$\||x|^k\nabla^k\bn\|_{L^\infty(D^2)}<+\infty$, we can apply 
lemma~\ref{lma-A.9}   in each annulus $D^2_{2^{-i}}\setminus D^2_{2^{-i-1}}$ to deduce that $\bL$ is smooth in $D^2\setminus\{0\}$ and that 
\be
\label{IV.10}
sup_{x\in D^2}|x|\ |\nabla \bL(x)|<+\infty\quad .
\ee     
Like in the previous section we introduce the cut-off function
$\chi$ equals to 1 on $D^2_{1/2}$ and compactly supported in $D^2$. We introduce $\bg_1$
and $\bg_2$ like in (\ref{III.c47}) and (\ref{III.c48}) and consider $\bv_1$ and $\bv_2$  the solution respectively of (\ref{III.7}) given by lemma~\ref{lm-III.1} for $\bg=\bg_1$ and  the solution of (\ref{III.32}) given by lemma~\ref{lm-III.3} for $\bg=\bg_2$. $\bv_1$ and $\bv_2$ satisfy (\ref{III.c51}) and 
(\ref{III.c52}). Therefore we have in particular that $\nabla \bv_1$ and $\nabla\bv_2$ are in
 $L^{2,\infty}$ and, like for $\bL$, since ${\mathcal L}_{\bn}\bv_i$ equals $0$ on $D^2_{1/2}$ 
 and since $\bg_1$ and $\bg_2$ are smooth, we have for $i=1,2$
 \be
 \label{IV.11}
sup_{x\in D^2}|x|\ |\nabla \bv_i(x)|<+\infty\quad .
\ee     
Denote $\bw:=\bH-\bv_1-\bv_2-\bL$. It is an $L^2$ solution to ${\mathcal L}_{\bn}\bw=0$ which
is smooth in $D^2\setminus\{0\}$ and equal to 0 on $\p D^2$. We claim that $\bw$ is identically
0 on $D^2$. 

For $r>0$ we denote $\chi_r(x)=\chi(x/r)$.
 Let $\bp_{i}$ be a sequence of normalised eigenvectors
of ${\mathcal {L}}_{\bn}$ in $W^{1,2}_0(D^2,{\R}^m)$ and forming an Hilbert orthonormal Basis of $L^2(D^2,{\R}^m)$
(the existence of such a Basis is obtained by combining the result of Lemma~\ref{lm-III.1} and Hilbert-Schmidt theorem).
Denote $\la_i$ the corresponding eigenvalues. From lemma~\ref{lm-III.1} again we know that $\la_i\ne 0$.
We have
\be
\label{IV.12}
\begin{array}{l}
\ds\int_{D^2}(1-\chi_r)\ \bw\cdot\bp_i=\frac{1}{\la_i}\int_{D^2}(1-\chi_r)\ \bw\cdot{\mathcal L}_{\bn}\bp_i\\[5mm]
\ds=\frac{1}{\la_i}\int_{D^2}\nabla \chi_r\ \bw\cdot \lf[\nabla\bp_i-3\pi_{\bn}(\nabla\bp_i)-\bp_i\wedge\nabla^\perp\bn\rg]\\[5mm]
\ds\ -\frac{1}{\la_i}\int_{D^2}(1-\chi_r)\ \nabla\bw\cdot \lf[\nabla\bp_i-3\pi_{\bn}(\nabla\bp_i)-\bp_i\wedge\nabla^\perp\bn\rg]\\[5mm]
\ds=\frac{1}{\la_i}\int_{D^2}(1-\chi_r)\ \bp_i\cdot {\mathcal L}_{\bn}\bw+\frac{1}{\la_i}\int_{D^2}\nabla \chi_r\ \bw\cdot \lf[\nabla\bp_i-3\pi_{\bn}(\nabla\bp_i)\rg]\\[5mm]
\ds\ -\frac{1}{\la_i}\int_{D^2}\nabla\chi_r\ \lf[\nabla\bw-3\pi_{\bn}(\nabla\bw)-\bw\wedge\nabla^\perp\bn\rg]\cdot\bp_i
\end{array}
\ee
Since in the sense of distribution ${\mathcal L}_{\bn}\bw=0$, we have
\be
\label{IV.13}
-\int_{D^2}\nabla\chi_r\ \lf[\nabla\bw-3\pi_{\bn}(\nabla\bw)-\bw\wedge\nabla^\perp\bn\rg]=
\lf<{\mathcal L}_{\bn}\bw,\chi_r\rg>_{{\mathcal D}',C^\infty_0}=0
\ee
Thus, in the last term of (\ref{IV.12}), we can substract to $\bp_i$ the vector $\bc_{r,i}$ which is the average of $\bp_i$ on $D^2_r\setminus D^2_{r/2}$, and we get
\be
\label{IV.14}
\begin{array}{l}
\ds\int_{D^2}(1-\chi_r)\ \bw\cdot\bp_i=\frac{1}{\la_i}\int_{D^2}\nabla \chi_r\ \bw\cdot \lf[\nabla\bp_i-3\pi_{\bn}(\nabla\bp_i)\rg]\\[5mm]
\ds\ -\frac{1}{\la_i}\int_{D^2}\nabla\chi_r\ \lf[\nabla\bw-3\pi_{\bn}(\nabla\bw)-\bw\wedge\nabla^\perp\bn\rg]\cdot(\bp_i-\bc_{r,i})\quad .
\end{array}
\ee
Denote $\nu(r)=sup_{r/2<|x|<r}|x|^2\ |\nabla\bw|(x)+|x|\ |w|(x)$. With this notation we control the
right-hand-side in the following way. In one hand, using Cauchy-Schwartz,
\be
\label{IV.15}
\begin{array}{l}
\ds\lf| \frac{1}{\la_i}\int_{D^2}\nabla \chi_r\ \bw\cdot \lf[\nabla\bp_i-3\pi_{\bn}(\nabla\bp_i)\rg]\rg|\\[5mm]
\ds\quad\le C_i\ \int_{D^2_r\setminus D^2_{r/2}}\frac{\nu(r)}{r^2}\ |\nabla\bp_i|\\[5mm]
\ds\quad\le C_i\ \nu(r)\ \lf[\int_{D^2_r\setminus D^2_{r/2}} \frac{|\nabla \bp_i|^2}{|x|^2}\rg]^\frac{1}{2}
\end{array}
\ee
Using lemma~\ref{lma-A.5}, and the fact that $\nu(r)$ converges  to zero as $r$ converges to zero
(This is obtained combining (\ref{IV.6}), (\ref{IV.10}) and (\ref{IV.11})) we have that the left-hand-side
of (\ref{IV.15}) converges to 0 as $r$ goes to 0.
In the other hand, we have, using Cauchy-Schwartz and Poincar\'e inequalities 
\be
\label{IV.16}
\begin{array}{l}
\ds\lf|-\frac{1}{\la_i}\int_{D^2}\nabla\chi_r\ \lf[\nabla\bw-3\pi_{\bn}(\nabla\bw)-\bw\wedge\nabla^\perp\bn\rg]\cdot(\bp_i-\bc_{r,i})\rg|\\[5mm]
\ds \quad\le C_i\ \int_{D^2_r\setminus D^2_{r/2}} \lf[\frac{\nu(r)}{r^3}+\frac{\nu(r)\ \delta(r)}{r^3}\rg]\
|\bp_i-\bc_{r,i}|\\[5mm]
\ds\quad\le C_i\ \lf[\nu(r)+\nu(r)\ \delta(r)\rg]\frac{1}{r^2}\ \lf[ \int_{D^2_r\setminus D^2_{r/2}} |\bp_i-\bc_{r,i}|^2\rg]^\frac{1}{2}\\[5mm]
\ds \quad\le C_i\ \lf[\nu(r)+\nu(r)\ \delta(r)\rg]\ \lf[ \int_{D^2_r\setminus D^2_{r/2}} \frac{|\nabla\bp_i|^2}{|x|^2}\rg]^\frac{1}{2}
\end{array}
\ee
Again using lemma~\ref{lma-A.5}, and the fact that $\nu(r)$ and $\delta(r)$ converge  to zero as $r$ converges to zero we obtain that the left-hand-side
of (\ref{IV.16}) also converges to 0 as $r$ goes to 0. So finally combining (\ref{IV.14}), (\ref{IV.15})
and (\ref{IV.16}) we have that $\int_{D^2}(1-\chi_r)\ \bw\cdot\bp_i$ converges to zero as $r$ goes to zero
which implies that $\int_{D^2}\bw\cdot\bp_i=0$. Since this holds for every $i$ and since $\bp_i$
realizes an orthonormal basis of $L^2$ we obtain that $\bw$ is identically zero and hence
\be
\label{IV.17}
\bH=\bL+\bv_1+\bv_2
\ee
Because of lemma~\ref{lma-A.9}, since ${\mathcal L}_{\bn}\bv_i$ equals to zero on $D^2_{1/2}$, we deduce that $\bv_1$ and $\bv_2$ are smooth on $D^2_{1/2}$. Thus it remains to study the assymptotic expansion of $\bL$ at the origin. First we observe that since $\nabla\bL\in L^{2,\infty}$, we have that
$\nabla \bH\in L^{2,\infty}$ and using (\ref{II.8a}) we deduce that 
\be
\label{IV.18}
\Delta\nabla\Phi=4e^{2\la}\bH\nabla\la+2e^{2\la}\nabla\bH\in \cap_{p<2} L^p\quad.
\ee
Since $e^\la=|\nabla\Phi|$ we deduce from (\ref{IV.18}) that $\nabla e^\la\in L^q$ for every $q<+\infty$.
Bootstraping this information again in (\ref{IV.18}) we deduce that $\Delta\nabla\bp\in L^{2,\infty}$
which implies in particular that $\nabla^2\bn$ is in $L^{2,\infty}$. From \cite{Hel} (proof of theorem 5.1.1)
we see that the Coulomb framing $(\bbe_1,\bbe_2)$ has the same regularity as $\bn$ : 
$\nabla^2\bbe_i\in L^{2,\infty}$. In particular this implies that $\bbe_i\in C^{0,\al}$ for every $0<\al<1$.
We claim now that $e_i(0)\cdot \bH_0=0$.
We have, using the fact that $\bH\cdot\bbe\equiv 0$,
\be
\label{IV.19}
\begin{array}{l}
\ds -4\pi\bbe_i(0)\cdot\bH_0\delta_0=\bbe_i\cdot div(\nabla \bH-3\pi_{\bn}(\nabla\bH)-\bH\wedge\nabla^\perp\bn)\\[5mm]
\ds\quad=div(\bbe_i\cdot\nabla\bH-3\bbe_i\cdot\bH\wedge\nabla^\perp\bn)-\nabla\bbe_i\cdot\lf[\nabla \bH-3\pi_{\bn}(\nabla\bH)-\bH\wedge\nabla^\perp\bn\rg]\\[5mm]
\ds=div(-\bH\cdot\nabla\bbe_i-3\bbe_i\cdot\bH\wedge\nabla^\perp\bn)-\nabla\bbe_i\cdot\lf[\nabla \bH-3\pi_{\bn}(\nabla\bH)-\bH\wedge\nabla^\perp\bn\rg]
\end{array}
\ee
Observe now that the right-hand-side of (\ref{IV.19}) is an $L^p$ function for some $p>1$ and that this function should be proportional to the Dirac mass at the origin. This implies that the coefficient $4\pi\bbe_i(0)\cdot\bH_0$ is zero. Let now $\bR:=\bL-\bH_0\ \log\, |x|$. We have
\be
\label{IV.20}
{\mathcal L}_{\bn}\bR=-3 div(\pi_{\bn}(\bH_0)\ \nabla\log|x|)-\nabla \log|x|\ \bH_0\wedge\nabla^\perp\bn
\ee
Since  $\pi_{\bn}(\bH_0)=(\bH_0\cdot\bbe_1)\,\bbe_1+(\bH_0\cdot\bbe_2)\,\bbe_2=
(\bH_0\cdot(\bbe_1-\bbe_1(0)))\,\bbe_1+(\bH_0\cdot(\bbe_2-\bbe_2(0)))\,\bbe_2$ and since
$\bbe_i$ are in $C^{0,\al}$ for every $\al<1$ we have then than $r^{-1}\pi_{\bn}(\bH_0)\in L^p$ for
every $p<+\infty$.  Thus we have proved that ${\mathcal L}_{\bn}\bR\in W^{-1,p}$ for every $p<+\infty$.
Arguing like in the proof of lemma~\ref{lm-III.1} we have that $\bR\in \cap_{p<+\infty}W^{1,p}$.
Hence we have proved that $\bH-\bH_0\ \log|x|$ is in $C^{0,\al}$ for every $\al<1$ and this finishes
the proof of theorem~\ref{th-III.1}.\hfill$\Box$

\appendix
\section{Appendix}
\reset
\begin{Lma}
\label{lm-III.1}
There exists $\ep_0>0$ such that for every $0<\ep<\ep_0$ , there exists a constant $C>0$ independent of $\ep$ such that the following holds. Let $\bn$ from $D^2$ into the space of
 unit $m-2-$vectors in ${\R}^m$ satisfying
\be
\label{III.6z}
\int_{D^2}|\nabla \bn|^2\ dx\le\ep\quad.
\ee
Let $\bg$ be an arbitrary distribution in the Sobolev space $H^{-1}(D^2,{\R}^m)$ dual to $W^{1,2}_0(D^2,{\R}^m)$, then there exists a unique
map $\bv$ in $W^{1,2}_0(D^2,{\R}^m)$ satisfying
\be
\label{III.7z}
\lf\{
\begin{array}{l}
\ds \Delta\bv-3\ div(\pi_{\bn}(\nabla \bv))-div\lf(\bv\wedge\nabla^\perp\bn\rg) =\bg\quad\mbox{ in } D^2\\[5mm]
\ds \bv=0\quad\quad\mbox{ on }\p D^2\quad ,
\end{array}
\rg.
\ee 
and
\be
\label{III.8z}
\int_{D^2}|\nabla \bv|^2\le C\ \|\bg\|_{H^{-1}}^2\quad .
\ee
Moreover the operator ${\mathcal{L}}_{\bn}^{-1}$ which to $\bg$ assigns $\bv$ satisfying (\ref{III.7z}) is selfadjoint and  compact
from $L^2(D^2,{\R}^m)$ into itself. \hfill$\Box$
\end{Lma}

Before proving lemma~\ref{lm-III.1} we shall prove first the following.
\begin{Lma}
\label{lm-III.2}
There exists $\ep_0>0$ such that for every $0<\ep<\ep_0$ , there exists a constant $C>0$ independent of $\ep$
 such that the following holds. Let $\bn$ from $D^2$ into the space of
 unit $m-2-$vectors in ${\R}^m$ satisfying
\be
\label{III.6}
\int_{D^2}|\nabla \bn|^2\ dx\le\ep\quad.
\ee
Let $\bg$ be an arbitrary distribution in the Sobolev space $H^{-1}(D^2,{\R}^m)$ dual to $W^{1,2}_0(D^2,{\R}^m)$, then there exists a unique
map $\bv$ in $W^{1,2}_0(D^2,{\R}^m)$ satisfying
\be
\label{III.7}
\lf\{
\begin{array}{l}
\ds \Delta\bv-3\ div(\pi_{\bn}(\nabla \bv))=\bg\quad\mbox{ in } D^2\\[5mm]
\ds \bv=0\quad\quad\mbox{ on }\p D^2\quad ,
\end{array}
\rg.
\ee 
and
\be
\label{III.8}
\int_{D^2}|\nabla \bv|^2\le C\ \|\bg\|_{H^{-1}}^2\quad .
\ee
\hfill$\Box$
\end{Lma}
{\bf Proof of lemma~\ref{lm-III.2}.}

First of all we show that under the assumption (\ref{III.6}) of the lemma, the following implication holds for every
$\bC$ in $W^{1,2}(D^2,{\R}^m)$ 
\be
\label{III.9}
\lf\{
\begin{array}{l}
\ds\Delta\bC-3\ div(\pi_{\bn}(\nabla \bC))=0\quad\mbox{ in }D^2\\[5mm]
\ds\bC=0\quad\quad\mbox{ on }\p D^2
\end{array}
\Longrightarrow\quad\quad \bC\equiv 0
\rg.
\ee
The proof of (\ref{III.9}) goes as follows. Since 
\[
div(\nabla \bC-3\pi_{\bn}(\nabla\bC))=0\quad,
\]
 from Poincar\'e lemma, there exists
$\bD$ in $W^{1,2}(D^2,{\R}^m)$ satisfying
\be
\label{III.10}
\nabla^\perp \bD=\nabla \bC-3\pi_{\bn}(\nabla\bC)\quad.
\ee
This implies in particular that $\bD$ is a $W^{1,2}$ solution of the following equation
\be
\label{III.11}
\lf\{
\begin{array}{l}
\ds\Delta\bD=3\sum_{k=1}^m\nabla^\perp C^k\cdot\nabla(e^k_1\ \bbe_1)+3\sum_{k=1}^m\nabla^\perp C^k\cdot\nabla(e^k_2\ \bbe_2)\quad
\mbox{ in }D^2\\[5mm]
\ds \frac{\p \bD}{\p\nu}=0\quad\quad\mbox{ on }\p D^2\quad ,
\end{array}
\rg.
\ee
where $C^k$ are the coordinates of $\bC$ in the canonical basis of ${\R}^m$, $(\bbe_1,\bbe_2)$ is an orthonormal basis of the 2-dimensional
subspace defined by $\bn$ given by lemma 5.1.4 in \cite{Hel}. From this lemma we have in particular the existence of a constant $C$ such that
\be
\label{III.12}
\int_{D^2}|\nabla\bbe_1|^2+|\nabla\bbe_2|^2\ dx\le C\ \int_{D^2}|\nabla \bn|^2\ dx\quad .
\ee
Using now the Neuman boundary condition version of Wente's lemma given in \cite{Hel} (lemma 3.1.2) we have, 
\be
\label{III.13}
\begin{array}{rl}
\ds\int_{D^2} |\nabla \bD|^2\ dx &\ds\le C_1\ \lf[\int_{D^2}|\nabla\bbe_1|^2+|\nabla\bbe_2|^2\ dx\rg] \int_{D^2}|\nabla \bC|^2\ dx\\[5mm]
\ds &\ds\le C_1\ep\ \int_{D^2}|\nabla \bC|^2\ dx
\end{array}
\ee
(the boundary condition studied in \cite{Hel} lemma 3.1.2 is the Dirichlet condition but modulo a 
slight modification the $W^{1,2}$ estimate can be obtained also for the Neuman boundary condition
by classical arguments in elliptic theory)

Observe that (\ref{III.10}) implies that 
\be
\label{III.14}
|\nabla\bD|^2=| \pi_T(\nabla \bC)|^2+4|\pi_{\bn}(\nabla \bC)|^2\ge|\nabla\bC|^2\quad,
\ee
where $\pi_T$ denotes the orthogonal projection on the 2-plane in ${\R}^m$ defined by $\bn$.
Combining then (\ref{III.13}) and (\ref{III.14}) we obtain, for $\ep<1/(2C_1)$,
 that $\bC\equiv 0$ which proves the implication (\ref{III.9}).

Let now $\bg$ in $H^{-1}(D^2,{\R}^m)$ and let $\bB$ solving
\be
\label{III.15}
\lf\{
\begin{array}{l}
\ds\Delta\bB=\bg\quad\quad\mbox{ in }D^2\\[5mm]
\ds\bB=0\quad\quad\mbox{ on }\p D^2
\end{array}
\rg.
\ee
We claim that there exists $(\bA,\bF)$ solutions to 
\be
\label{III.16}
\lf\{
\begin{array}{l}
div\, \bF=div\lf(\pi_T(\nabla^\perp\bA)-\frac{1}{2}\pi_{\bn}(\nabla^\perp A)\rg)\quad\quad\mbox{ in } D^2\\[5mm]
curl\,\bF=-curl\lf(\pi_T(\nabla \bB)-\frac{1}{2}\pi_{\bn}(\nabla \bB)\rg)\quad\quad\mbox{ in } D^2\\[5mm]
\bF\cdot\nu=0\quad\quad\mbox{ on }\p D^2
\end{array}
\rg.
\ee
where $\bA$ is the $curl-$part in the  Hodge decomposition of $\pi_T(\bF)-2\pi_{\bn}(\bF)$ given by :
\be
\label{III.17}
\lf\{
\begin{array}{l}
\ds -\Delta \bA=curl\lf(\pi_T(\bF)-2\pi_{\bn}(\bF)\rg)\quad\quad\mbox{ in }D^2\\[5mm]
\ds \bA=0\quad\quad\mbox{ on }\p D^2
\end{array}
\rg.
\ee
(remark that $\bF(x)$ is an element of ${\R}^2\otimes{\R}^m$ and that $\bF\cdot\nu\in{\R}^m$ in (\ref{III.16}) has to do with
the scalar product between the unit exterior normal to $\p D^2$ and  the ${\R}^2$ part in $\bF$).
 The existence of a solution $(\bA,\bF)$ of the system (\ref{III.16}), (\ref{III.17}) is again a consequence of Wente's estimate :
we write in one hand
\be
\label{III.18}
div\lf(\pi_T(\nabla^\perp\bA)-\frac{1}{2}\pi_{\bn}(\nabla^\perp A)\rg)=\frac{3}{2}\sum_{k=1}^m\nabla^\perp A^k\cdot\nabla(e_1^k\ \bbe_1)+
\frac{3}{2}\sum_{k=1}^m\nabla^\perp A^k\cdot\nabla(e_2^k\ \bbe_2)
\ee
and in the other hand
\be
\label{III.19}
curl\lf(\pi_T(\nabla \bB)-\frac{1}{2}\pi_{\bn}(\nabla \bB)\rg)=\frac{3}{2}\sum_{k=1}^m\nabla B^k\cdot\nabla^\perp(e_1^k\ \bbe_1)+
\frac{3}{2}\sum_{k=1}^m\nabla B^k\cdot\nabla^\perp(e_2^k\ \bbe_2)
\ee
where $A^k$ and $B^k$ are the coordinates of respectively $\bA$ and $\bB$. Therefore, using Wente's estimate, we have the following
a-priori inequalities
\be
\label{III.20}
\begin{array}{rl}
\ds\int_{D^2}|\bF|^2&\ds\le C_2\ \lf[\int_{D^2}|\nabla\bbe_1|^2+|\nabla\bbe_2|^2 dx\rg]\ \int_{D^2}|\nabla \bA|^2+|\nabla \bB|^2 dx\\[5mm]
\ds &\ds\le C_2\ \ep\ \int_{D^2}|\nabla \bA|^2+|\nabla \bB|^2 dx\quad.
\end{array}
\ee
From (\ref{III.17}) and standard elliptic estimates we have
\be
\label{III.21} 
\int_{D^2}|\nabla \bA|^2\ dx\le C_3 \int_{D^2}|\bF|^2\ dx\quad .
\ee
Thus for $C_3\,C_2\,\ep<1/2$, a standard fixed point argument gives the existence and uniqueness of $(\bA,\bF)$ satisfying (\ref{III.16})
and (\ref{III.17}).
Since 
\[
div\lf(\, \bF-\lf(\pi_T(\nabla^\perp\bA)-\frac{1}{2}\pi_{\bn}(\nabla^\perp \bA)\rg)\rg)=0\quad,
\]
there exists $\bC$ in $W^{1,2}_0(D^2,{\R}^m)$
satisfying
\be
\label{III.22}
\bF-\lf(\pi_T(\nabla^\perp\bA)-\frac{1}{2}\pi_{\bn}(\nabla^\perp A)\rg)=\nabla^\perp \bC\quad .
\ee
 Then we deduce that
\be
\label{III.23}
\pi_T(\bF)-2\pi_{\bn}(\bF)=\nabla^\perp\bA+ \pi_T(\nabla^\perp\bC)-2\pi_{\bn}(\nabla^\perp\bC)\quad.
\ee
Taking the curl of this identity together with (\ref{III.17}) implies that $\bC$ solves the following
equation
\be
\label{III.24}
\lf\{
\begin{array}{l}
\ds\Delta\bC-3\ div(\pi_{\bn}(\nabla \bC))=-curl\lf(\pi_T(\nabla^\perp\bC)-2\pi_{\bn}(\nabla^\perp\bC)\rg)=0\quad\mbox{ in }D^2\\[5mm]
\ds\bC=0\quad\quad\mbox{ on }\p D^2\quad .
\end{array}
\rg.
\ee
(\ref{III.9}) implies then that $\bC\equiv 0$ therefore 
\be
\label{III.25}
\pi_T(\bF)-2\pi_{\bn}(\bF)=\nabla^\perp\bA\quad .
\ee
From (\ref{III.16}) there exists $\bv$ in $W^{1,2}_0({D^2},{\R}^m)$ satisfying
\be
\label{III.26}
\bF=-\pi_T(\nabla \bB)+\frac{1}{2}\pi_{\bn}(\nabla \bB)+\nabla\bv
\ee
Combining (\ref{III.25}) and (\ref{III.26}) we obtain that
\be
\label{III.27}
\nabla^\perp \bA=-\nabla \bB+\pi_T(\nabla\bv)-2\pi_{\bn}(\nabla\bv)\quad.
\ee
Comparing this identity with (\ref{III.15}) we deduce that $\bv$ solves (\ref{III.7})
and from (\ref{III.9}) we know that this is the unique solution. (\ref{III.8}) follows from
(\ref{III.20}) and (\ref{III.21}). This completes then the proof of lemma~\ref{lm-III.2}. 
\hfill$\Box$

\medskip

\noindent {\bf Proof of Lemma~\ref{lm-III.1}.}

Denote $\Delta^{-1}_0$ the continuous isomorphism from $H^{-1}(D^2,{\R}^m)$ into $W^{1,2}_0(D^2,{\R}^m)$ which to
a distribution $\bg$ in $H^{-1}(D^2,{\R}^m)$ assigns the solution $\bv$ of
\be
\label{III.28}
\lf\{
\begin{array}{l}
\Delta\,\bv=\bg\quad\quad\mbox{ in }D^2\\[5mm]
\bv=0\quad\quad\mbox{ on }\p D^2\quad .
\end{array}
\rg.
\ee
We have seen in lemma~\ref{lm-III.2} that the operator ${\mathcal{A}}_{\bn}\bv:=\Delta\bv-3\ div(\pi_{\bn}(\nabla \bn))$
is a continuous isomorphism from $W^{1,2}_0(D^2,{\R}^m)$ into $H^{-1}(D^2,{\R}^m)$ and the norm of ${\mathcal{A}}_{\bn}$
and ${\mathcal{A}}_{\bn}^{-1}$ are independent of $\bn$ satisfying (\ref{III.6}) for $\ep<\ep_0$ where $\ep_0$ is universal
given by lemma~\ref{lm-III.2}.
Our aim is to show the invertibility of the operator
\[
\Delta_0^{-1}{\mathcal{A}}_{\bn}(\bv)-\Delta_0^{-1}div(\bv\wedge\nabla^\perp\bn)
\]
from $W^{1,2}_0(D^2,{\R}^m)$ into itself, with norms independant of $\bn$ satisfying (\ref{III.6}) for $\ep<\ep_0$.
This is a direct consequence of the invertibility of $\Delta_0^{-1}{\mathcal{A}}_{\bn}$ from $W^{1,2}_0(D^2,{\R}^m)$ into itself
and the fact that by the mean of Wente estimate ( theorem 3.1.2 of \cite{Hel}) the operator $\Delta_0^{-1}div(\bv\wedge\nabla^\perp\bn)$
satisfies for every $\bv$ in $W^{1,2}_0(D^2,{\R}^m)$
\be
\label{III.29}
\|\Delta_0^{-1}div(\bv\wedge\nabla^\perp\bn)\|^2_{W^{1,2}}\le C\ \int_{D^2}|\nabla \bv|^2dx\ \int_{D^2}|\nabla \bn|^2dx\le C\ \ep\ 
\|\bv\|^2_{W^{1,2}}
\ee
We have then proved the first statement of lemma~\ref{lm-III.1} and it remains to show the compactness and the self-adjointness of 
${\mathcal{L}}_{\bn}^{-1}$ from $L^2$ into itself.
Compactness is clear since ${\mathcal{L}}_{\bn}^{-1}$ sends $H^{-1}(D^2,{\R}^m)$ into $W^{1,2}_0$ which is a compact subspace of $L^2$.
Let $\bg$ and $\bbh$ in $L^2(D^2,{\R}^m)$ (we can take them first smooth). Denote $\bv:={\mathcal{L}}_{\bn}^{-1}(\bg)$ and
 $\bw:={\mathcal{L}}_{\bn}^{-1}(\bbh)$ we have
\be
\label{III.30}
\begin{array}{l}
\ds\int_{D^2}\bg \cdot {\mathcal{L}}_{\bn}^{-1}(\bbh)\ dx\\[5mm]
\ds=\int_{D^2} \Delta\bv\cdot\bw-3\ div(\pi_{\bn}(\nabla \bv))\cdot\bw
-div\lf(\bv\wedge\nabla^\perp\bn\rg)\cdot\bw\ dx\\[5mm]
 \ds=\int_{D^2} \bv\cdot\Delta\bw+3\pi_{\bn}(\nabla \bv)\cdot\nabla\bw+(\bv\wedge\nabla^\perp\bn)\cdot\nabla\bw\ dx\\[5mm]
\ds=\int_{D^2} \bv\cdot\Delta\bw+3\nabla \bv\cdot\pi_{\bn}(\nabla\bw)-\bv\cdot(\nabla\bw\wedge\nabla^\perp\bn)\\[5mm]
\ds=\int_{D^2} \bv\cdot\Delta\bw-3\bv\cdot div\lf(\pi_{\bn}(\nabla\bw)\rg)-\bv\cdot div\lf(\bw\wedge\nabla^\perp\bn\rg)\\[5mm]
\ds=\int_{D^2}\bv\cdot\bbh=\int_{D^2}{\mathcal{L}}_{\bn}^{-1}(\bg) \cdot \bbh\quad.
\end{array}
\ee
This shows the self-adjointness of ${\mathcal{L}}_{\bn}^{-1}$ and lemma~\ref{lm-III.1} is proved.
\hfill$\Box$

\medskip

We now extend the two previous lemma to $L^1$ datas. We have first

\begin{Lma}
\label{lm-III.3}
There exists $\ep_0>0$ such that for every $0<\ep<\ep_0$ , there exists a constant $C>0$ independent of $\ep$ such that the following holds. Let $\bn$ from $D^2$ into the space of
 unit $m-2-$vectors in ${\R}^m$ satisfying
\be
\label{III.31}
\int_{D^2}|\nabla \bn|^2\ dx\le\ep\quad.
\ee
Let $\bg$ be an arbitrary map in $L^1(D^2,{\R}^m)$, then there exists a unique
map $\bv$ with $\nabla\bv$ in $L^{2,\infty}(D^2,{\R}^2\otimes{\R}^m)$  satisfying
\be
\label{III.32}
\lf\{
\begin{array}{l}
\ds \Delta\bv-3\ div(\pi_{\bn}(\nabla \bv))-div\lf(\bv\wedge\nabla^\perp\bn\rg) =\bg\quad\mbox{ in } D^2\\[5mm]
\ds \bv=0\quad\quad\mbox{ on }\p D^2\quad ,
\end{array}
\rg.
\ee 
and
\be
\label{III.33}
\|\nabla\bv\|_{L^{2,\infty}(D^2)}\le C\ \|\bg\|_{L^1(D^2)}\quad .
\ee
\hfill$\Box$
\end{Lma}

Before proving lemma~\ref{lm-III.3} we shall prove first the following.
\begin{Lma}
\label{lm-III.4}
There exists $\ep_0>0$ such that for every $0<\ep<\ep_0$ , there exists a constant $C>0$ independent of $\ep$
 such that the following holds. Let $\bn$ from $D^2$ into the space of
 unit $m-2-$vectors in ${\R}^m$ satisfying
\be
\label{III.34}
\int_{D^2}|\nabla \bn|^2\ dx\le\ep\quad.
\ee
Let $\bg$ be an arbitrary map in $L^1(D^2,{\R}^m)$, then there exists a unique
map $\bv$ with $\nabla\bv$ in $L^{2,\infty}(D^2,{\R}^2\otimes{\R}^m)$  satisfying
\be
\label{III.35}
\lf\{
\begin{array}{l}
\ds \Delta\bv-3\ div(\pi_{\bn}(\nabla \bv))=\bg\quad\mbox{ in } D^2\\[5mm]
\ds \bv=0\quad\quad\mbox{ on }\p D^2\quad ,
\end{array}
\rg.
\ee 
and
\be
\label{III.36}
\|\nabla\bv\|_{L^{2,\infty}(D^2)}\le C\ \|\bg\|_{L^1(D^2)}\quad .
\ee
\hfill$\Box$
\end{Lma}
{\bf Proof of lemma~\ref{lm-III.4}.}
We first recall the following generalization of Wente's estimate established in \cite{Hel}
(theorem 3.4.5).
Let $a$ and $b$ be two functions on $D^2$ with $\nabla a\in L^{2,\infty}(D^2)$ and
 $\nabla b\in L^2(D^2)$, then there is a unique solution $\varphi$ in $W^{1,2}_0(D^2,{\R})$
 satisfying
 \be
 \label{III.37}
 \lf\{
 \begin{array}{l}
 \ds\Delta\varphi=\nabla a\cdot\nabla^\perp b\quad\quad\mbox{ in }D^2\\[5mm]
 \varphi=0\quad\quad\mbox{ on }\p D^2\quad,
 \end{array}
 \rg.
 \ee
 and the following estimate holds : there exists a positive constant $C$ independent of $a$
 and $b$ such that
 \be
 \label{III.38}
 \|\nabla\varphi\|_{L^2(D^2)}\le C\|\nabla a\|_{L^{2,\infty}(D^2)}\ \|\nabla b\|_{L^2(D^2)}\quad.
 \ee
 The proof of lemma~\ref{lm-III.4} follows then step by step the proof of lemma~\ref{lm-III.2}
 by replacing the original Wente estimate by the estimate (\ref{III.38}) and the $L^2$ norms
 of $\nabla\bC$, $\nabla\bB$, $\nabla\bA$, $\bF$ and $\nabla\bv$ replaced by their $L^{2,\infty}$ norms.
\hfill$\Box$
 
 \medskip
 
\noindent{\bf Proof of lemma~\ref{lm-III.3}.}
 
 Let $\bg$ be in $L^1(D^2,{\R}^m)$. Denote again ${\mathcal{A}}_{\bn}\bv:=\Delta\bv-3\ div(\pi_{\bn}(\nabla \bv))$. Using 
 lemma~\ref{lm-III.4} we first get the existence of $\bv_0$ with $\nabla\bv_0\in L^{2,\infty}$ satisfying
 \be
 \label{III.39}
 \lf\{
 \begin{array}{l}
 {\mathcal{A}}_{\bn}\bv_0=\bg\quad\quad\mbox{ in }D^2\\[5mm]
 \bv_0=0\quad\quad\mbox{ on }\p D^2
 \end{array}
 \rg.
 \ee
 Then we construct by induction the following sequence $\bv_k$ with $\bv_0$ given by
 (\ref{III.39}) and $\bv_k$ for $k\ge 1$ will be the element of $W^{1,2}_0(D^2,{\R}^m)$ solving
 \be
 \label{III.40}
 \Delta_0^{-1}{\mathcal{A}}_{\bn}(\bv_k)=\Delta_0^{-1}div(\bv_{k-1}\wedge \nabla^\perp \bn)\quad\quad,
 \ee
 where $\Delta_0^{-1}$ is the operator introduced in the proof of lemma~\ref{lm-III.1} by (\ref{III.28}).
 This problem admits indeed a solution for the following reason : let $\bv_{k-1}$ being 
 given, with the assumption that $\nabla \bv_{k-1}$ in $L^{2,\infty}$, using (\ref{III.38}) we have
 \be
 \label{III.41}
 \|\Delta_0^{-1}div(\bv_{k-1}\wedge \nabla^\perp \bn)\|_{W^{1,2}}\le C\ \|\nabla \bv_{k-1}\|_{L^{2,\infty}}
 \|\nabla\bn\|_{L^2}\quad.
 \ee
Moreover we have seen in the proof of lemma~\ref{lm-III.1}, using
lemma~\ref{lm-III.2}, that $ \Delta_0^{-1}{\mathcal{A}}_{\bn}$ is a continuous isomorphism of $W^{1,2}_0$.
Combining (\ref{III.41}) and this last fact we get the existence and uniqueness of $\bv_{k}$, once
$\bv_{k-1}$ is known, and the following inequality holds
\be
\label{III.42}
\|\nabla\bv_k\|_{L^{2,\infty}}\le \|\nabla\bv_k\|_{L^2}\le C\ \|\nabla \bv_{k-1}\|_{L^{2,\infty}}
 \|\nabla\bn\|_{L^2}\quad.
 \ee
Thus, under the assumption that   $C \|\nabla\bn\|_{L^2}<1/2$ the serie $\sum_{k=0}^n\bv_k$ converges
to a limit $\bv=\sum_{k=0}^\infty\bv_k$ solving (\ref{III.32}) and (\ref{III.33}) holds. The uniqueness of
$\bv$ follows from the uniqueness of the solution of ${\mathcal{A}}_{\bn}(\bv)=\bg$ for arbitrary $\bg\in L^1$ established in lemma~\ref{lm-III.4} and the previous considerations with
$C \|\nabla\bn\|_{L^2}$ being small enough. \hfill$\Box$
\begin{Lma}
\label{lma-A.5}
There exists $\ep_0>0$ such that for every $0<\ep<\ep_0$ , there exists a constant $C>0$ independent of $\ep$
 such that the following holds. Let $\bn$ from $D^2$ into the space of
 unit $m-2-$vectors in ${\R}^m$ satisfying
\be
\label{III.43}
\||x|\,|\nabla\bn|(x)\|_{L^\infty(D^2)}+\int_{D^2}|\nabla \bn|^2\ dx\le\ep\quad.
\ee
Let $\bg$ be an arbitrary function in $L^2(D^2,{\R}^m)$ , let $\bv$ be the unique
map in $W^{1,2}_0(D^2,{\R}^m)$ given by lemma~\ref{lm-III.2} satisfying
\be
\label{III.44}
\lf\{
\begin{array}{l}
\ds \Delta\bv-3\ div(\pi_{\bn}(\nabla \bv))-div\lf(\bv\wedge\nabla^\perp\bn\rg)=\bg\quad\mbox{ in } D^2\\[5mm]
\ds \bv=0\quad\quad\mbox{ on }\p D^2\quad ,
\end{array}
\rg.
\ee 
and denote $\bv_0$ and $\bv_\perp$  the  maps from $D^2$ into ${\R}^m$ given by
\be
\label{III.45}
\bv_0(x)=\frac{1}{2\pi |x|}\int_{\p B_{|x|}(0)} \bv\quad\mbox{ and }\quad\bv_\perp=\bv-\bv_0\quad.
\ee
Then the following inequality holds
\be
\label{III.46}
\|\nabla \bv_0\|^2_{L^\infty(D^2)} + \int_{D^2}\frac{|\nabla\bv_\perp|^2}{|x|^2}+ \int_{D^2}{|\nabla^2\bv_\perp|^2
}\le C\ \int_{D^2}|\bg|^2
\ee
\end{Lma}
Before proving lemma~\ref{lma-A.5} we need to establish the following intermediate result.
\begin{Lma}
\label{lma-A.6}
Let $a$ and $b$ be two functions respectively in $W^{2,2}(D^2,{\R})$ and $W^{1,2}(D^2,{\R})$
 such that $b\in C^1(D^2\setminus\{0\})$
and
\[
\sup_{x\in D^2\setminus\{0\}}|x|\ |\nabla b|(x)<+\infty
\]
Let $\vp$ be the solution in $W^{1,2}$ of 
\be
\label{III.47}
\lf\{
\begin{array}{l}
\ds\Delta \varphi=\frac{\p a}{\p x_1}\,\frac{\p b}{\p x_2}-\frac{\p a}{\p x_2}\,\frac{\p b}{\p x_1}\quad\mbox{ in }D^2\\[5mm]
\ds \varphi=0\quad\quad\mbox{ on }\p D^2
\end{array}
\rg.
\ee
Denote $\varphi_0$ and $\varphi_\perp$ the functions on $D^2$ given by
\be
\label{III.48}
\varphi_0(x)=\frac{1}{2\pi |x|}\int_{\p B_{|x|}(0)} \varphi \quad\mbox{ and }\quad \varphi_\perp=\varphi-\vp_0\quad
\ee
Then the following inequality holds
\be
\label{III.49}
\begin{array}{l}
\ds\|\nabla \vp_0\|_{L^\infty(D^2)}^2+\int_{D^2}\frac{|\nabla \vp_\perp|^2}{|x|^2}+\int_{D^2}|\nabla^2 \vp_\perp|^2\\[5mm]
\ds\le C\ \lf[\||x|\,|\nabla b|(x)\|^2_\infty+\int_{D^2}|\nabla b|^2\rg]\ \lf[\int_{D^2}\frac{|\nabla a_\perp|^2}{|x|^2}
+\|\nabla a_0\|^2_{\infty}\rg]\quad.
\end{array}
\ee
\end{Lma}
{\bf Proof of lemma~\ref{lma-A.6}.}
Since $\vp_0$ is the first term in the Fourier decomposition of $\vp$ (for the angle variable), we have
that
\be
\label{III.50}
\begin{array}{l}
\ds\Delta \vp_0=\frac{\p a_0}{\p x_1}\,\frac{\p b_0}{\p x_2}-\frac{\p a_0}{\p x_2}\,\frac{\p b_0}{\p x_1}\\[5mm]
\ds\quad+\lf(\frac{\p a_\perp}{\p x_1}\,\frac{\p b_\perp}{\p x_2}-\frac{\p a_\perp}{\p x_2}\,\frac{\p b_\perp}{\p x_1}\rg)_0
\end{array}
\ee
Indeed it is clear that $\frac{\p a_\perp}{\p x_1}\,\frac{\p b_0}{\p x_2}-\frac{\p a_\perp}{\p x_2}\,\frac{\p b_0}{\p x_1}$
as well as $\frac{\p a_0}{\p x_1}\,\frac{\p b_\perp}{\p x_2}-\frac{\p a_0}{\p x_2}\,\frac{\p b_\perp}{\p x_1}$ have no $0-$component.
This comes from the fact that $a_0$ and $b_0$ depends only on $|x|$ which imply
the identities
\be
\label{III.51}
\begin{array}{l}
\ds\frac{\p a_\perp}{\p x_1}\,\frac{\p b_0}{\p x_2}-\frac{\p a_\perp}{\p x_2}\,\frac{\p b_0}{\p x_1}=\frac{1}{r}\frac{\p a_\perp}{\p \theta}\  \dot{b}_0(r)\\[5mm]
\ds\frac{\p a_0}{\p x_1}\,\frac{\p b_\perp}{\p x_2}-\frac{\p a_0}{\p x_2}\,\frac{\p b_\perp}{\p x_1}=\dot{a}_0(r)\
\frac{1}{r}\frac{\p b_\perp}{\p \theta}\quad .
\end{array}
\ee
It is also clear that, since both $a_0$ and $b_0$ only depends on $|x|$ that the first jacobien in 
(\ref{III.50}) is zero. Thus we have, using one more time (\ref{III.51}),
\be
\label{III.52}
\begin{array}{rl}
\ds\ddot{\vp}_0+\frac{\dot{\vp}}{r}=\Delta \vp_0 &\ds=\lf(\frac{\p a_\perp}{\p x_1}\,\frac{\p b_\perp}{\p x_2}-\frac{\p a_\perp}{\p x_2}\,\frac{\p b_\perp}{\p x_1}\rg)_0\\[5mm]
 &\ds=\lf(\frac{\p a_\perp}{\p x_1}\,\frac{\p b}{\p x_2}-\frac{\p a_\perp}{\p x_2}\,\frac{\p b}{\p x_1}\rg)_0
\end{array}
\ee
Denoting $h(r)$ the right-hand side of (\ref{III.52}), we then have that $\dot{\vp}_0=r^{-1}\int_0^rh(s)\ s\,ds$. Therefore $|\dot{\vp}_0|$
can be bounded in the following way :
\be
\label{III.53}
\begin{array}{rl}
\ds|\dot{\vp}_0|&\ds\le \frac{1}{2\pi r}\ \int_{B_r}|\nabla a_\perp|\,|\nabla b|\\[5mm]
&\ds \le C \||x|\,|\nabla b|(x)\|_\infty \lf[\int_{B_r}\frac{|\nabla a_\perp|^2}{|x|^2}\rg]^\frac{1}{2}\quad .
\end{array}
\ee
Regarding $\vp_\perp$ we have :
\be
\label{III.54}
\begin{array}{l}
\ds\Delta \vp_\perp =\lf(\frac{\p a}{\p x_1}\,\frac{\p b}{\p x_2}-\frac{\p a}{\p x_2}\,\frac{\p b}{\p x_1}\rg)_\perp\\[5mm]
 \ds\quad=\lf(\frac{\p a_0}{\p x_1}\,\frac{\p b}{\p x_2}-\frac{\p a_0}{\p x_2}\,\frac{\p b}{\p x_1}\rg)_\perp+\lf(\frac{\p a_\perp}{\p x_1}\,\frac{\p b}{\p x_2}-\frac{\p a_\perp}{\p x_2}\,\frac{\p b}{\p x_1}\rg)_\perp
\end{array}
\ee
In one hand we have
\be
\label{III.55}
\begin{array}{l}
\ds\int_{D^2}\lf|\lf(\frac{\p a_0}{\p x_1}\,\frac{\p b}{\p x_2}-\frac{\p a_0}{\p x_2}\,\frac{\p b}{\p x_1}\rg)_\perp\rg|^2\\[5mm]
\ds\le\int_{D^2}\lf|\frac{\p a_0}{\p x_1}\,\frac{\p b}{\p x_2}-\frac{\p a_0}{\p x_2}\,\frac{\p b}{\p x_1}\rg|^2\le C\|\nabla a_0\|^2_\infty\ \int_{D^2}|\nabla b|^2\quad.
\end{array}
\ee
In the other hand we have
\be
\label{III.56}
\begin{array}{l}
\ds\int_{D^2}\lf|\lf(\frac{\p a_\perp}{\p x_1}\,\frac{\p b}{\p x_2}-\frac{\p a_\perp}{\p x_2}\,\frac{\p b}{\p x_1}\rg)_\perp\rg|^2\\[5mm]
\ds\le\int_{D^2}\lf|\frac{\p a_\perp}{\p x_1}\,\frac{\p b}{\p x_2}-\frac{\p a_\perp}{\p x_2}\,\frac{\p b}{\p x_1}\rg|^2\le C\||x|\,|\nabla b|(x)\|_\infty^2\ \int_{D^2}\frac{|\nabla a_\perp|^2}{|x|^2}
\end{array}
\ee
Combining (\ref{III.53}), (\ref{III.55}) and (\ref{III.56}) we get (\ref{III.49}) and lemma~\ref{lma-A.6} is proved.\hfill$\Box$

\medskip

\noindent{\bf Proof of lemma~\ref{lma-A.5}.}
Let $\bA$, $\bB$ and $\bC$ be the solutions respectively of
\be
\label{III.57}
\lf\{
\begin{array}{l}
\ds\Delta\bA=\bg\quad\quad\mbox{ in }D^2\\[5mm]
\ds\bA=0\quad\quad\mbox{ on }\p D^2
\end{array}
\rg.
\ee
\be
\label{III.58}
\lf\{
\begin{array}{l}
\ds\Delta\bB=div(\bv\wedge\nabla^\perp\bn)\quad\quad\mbox{ in }D^2\\[5mm]
\ds\bB=0\quad\quad\mbox{ on }\p D^2
\end{array}
\rg.
\ee
and
\be
\label{III.59}
\lf\{
\begin{array}{l}
\ds\Delta\bC=3\,div(\pi_{\bn}(\nabla^\perp\bv))\quad\quad\mbox{ in }D^2\\[5mm]
\ds\frac{\p\bC}{\p\nu}=0\quad\quad\mbox{ on }\p D^2
\end{array}
\rg.
\ee
It is clear that $\bv=\bA+\bB+\bC$. Using standard elliptic estimates for equation (\ref{III.57}) and applying lemma~\ref{lma-A.6}  to equation
(\ref{III.58})  and the Neuman-boundary condition version of lemma~\ref{lma-A.6} to equation (\ref{III.59})
we obtain 
\be
\label{III.60}
\|\nabla \bA_0\|_\infty^2+\int_{D^2}\frac{|\nabla \bA_\perp|^2}{|x|^2}+\int_{D^2}|\nabla^2 \bA|^2\le C\ \int_{D^2}|\bg|^2\quad,
\ee
\be
\label{III.61}
\|\nabla \bB_0\|_\infty^2+\int_{D^2}\frac{|\nabla \bB_\perp|^2}{|x|^2}+\int_{D^2}|\nabla^2 \bB_\perp|^2\le C\ep\ \lf[\int_{D^2}\frac{|\nabla \bv_\perp|^2}{|x|^2}+
\|\nabla \bv_0\|_\infty^2\rg]
\ee
and
\be
\label{III.62}
\|\nabla \bC_0\|_\infty^2+\int_{D^2}\frac{|\nabla \bC_\perp|^2}{|x|^2}+\int_{D^2}|\nabla^2 \bC_\perp|^2\le C\ep\ \lf[\int_{D^2}\frac{|\nabla \bv_\perp|^2}{|x|^2}+
\|\nabla \bv_0\|_\infty^2\rg]
\ee
where we have used the fact (like in (\ref{III.11})) that $div(\pi_{\bn}(\nabla^\perp\bv))$ is a jacobian of the
 form $-\sum_{k,i}\nabla^\perp\bv^k\cdot\nabla(e^k_i\, \bbe_i)$
and $(\bbe_1,\bbe_2)$ is an orthonormal frame generating the 2-plane $\bn$ and given by lemma 5.1.4 of \cite{Hel}.
Combining (\ref{III.60}), (\ref{III.61}), (\ref{III.62}) and choosing $\ep$ small enough, we get (\ref{III.49}) and lemma~\ref{lma-A.5} is proved.\hfill$\Box$
\begin{Lma}
\label{lma-A.7}
There exists $\ep_0>0$ such that for every $0<\ep<\ep_0$ the following holds. 
Let $\bn$ be a $W^{1,2}$ map from $D^2$ into the space of unit $m-2-$vectors in ${\R}^m$ satisfying
\be
\label{III.63}
\int_{D^2}|\nabla \bn|^2\ dx\le\ep\quad.
\ee
Let $\bp$ be a $W^{1,2}$ eigenvector of ${\mathcal{L}}_{\bn}$ i.e. solution to the equation
\be
\label{III.64}
\lf\{
\begin{array}{l}
\ds\Delta\bp-3\ div(\pi_{\bn}(\nabla \bp))-div\lf(\bp\wedge\nabla^\perp\bn\rg)=\la\,\bp\quad\quad\mbox{ in } D^2
\\[5mm]
\ds\bp=0\quad\quad\mbox{on }\p D^2\quad,
\end{array}
\rg.
\ee
for some constant $\la\in{\R}$. Then, assuming furthermore that the gradient of $\bn$ is in the Lorentz space $L^{2,1}$, we have that $\bp$ is Lipschitz with second derivatives in $L^{2,1}$.\hfill$\Box$
\end{Lma}
\noindent{\bf Proof of lemma~\ref{lma-A.7}.}
First we prove that $\bp$ is in $W^{1,p}(D^2)$ for every $1\le p<+\infty$. Let $2<p<+\infty$ and denote $q$ the constant in $(1,2)$ given by $1/p=1/q-1/2$ in such a way that $W^{1,q}_0(D^2)$ embedds in $L^p$.
We denote $\bg:=\la\bp\in L^q(D^2)$ and we follow step by step
the proof of lemma~\ref{lm-III.2}  starting from \ref{III.15} replacing the assumption that $\bg$ is in $H^{-1}$ by the hypothesis 
that $\bg$ is in $L^q$. Precisely we first observe, using  classical elliptic estimate that
\be
\label{III.65}
\|\nabla \bB\|_{L^p(D^2)}\le C\ \|\bg\|_{L^q(D^2)}
\ee
Replacing Wente's inequalities by classical $L^q$ bound for Calderon-Zygmund operators we obtain  the following a-priori estimate 
\be
\label{III.66}
\begin{array}{rl}
\ds\|\bF\|_{L^p(D^2)}&\ds\le C \|\bF\|_{W^{1,q}}\le C\ \lf[\int_{D^2}|\nabla \bA|^q\ |\nabla \bbe|^q+|\nabla \bB|^q\ |\nabla \bbe|^q\rg]^\frac{1}{q}\\[5mm]
 &\ds\le C\ \|\nabla \bA\|_p\ \|\nabla\bbe\|_2+\|\nabla \bB\|\ \|\nabla\bbe\|_2\\[5mm]
  &\ds\le C\ep\ \lf[\|\nabla \bA\|_p+\|\nabla \bB\|_p\rg]\quad.
 \end{array}
 \ee
From this a-priori estimate, like in the proof of lemma~\ref{lm-III.2}, we get the existence of $\bv$ solving
${\mathcal{L}}_{\bn}\bv=\bg$ which this time in $W^{1,p}_0$ for any $p>2$. Because of the uniqueness 
given by lemma~\ref{lm-III.2} they all coincide with $\bp$ and we have then proved that $\bp\in\cap_{p<+\infty} W^{1,p}_0(D^2)$.

Let $x_0$ be a point in the interior of $D^2$. For any $\ep>0$ we can find a radius $\rho>0$ such that 
the $L^{2,1}$ Lorentz norm of $\nabla\bn$ on $B_\rho(x_0)$ is less than $\ep$
\be
\label{III.67}
\|\nabla \bn\|_{L^{2,1}(B_\rho(x_0))}\le\ep\quad\quad.
\ee
Consider a smooth cut-off function $\chi$ equals to 1 on $B_{1/2}(0)$ and equals to 0
outside $B_1(0)=D^2$. Denote $\bw$ the function on the two dimensional disk equals to
\[
\bw(x)=\bp(\rho x+x_0)\ \chi(x)\quad.
\]
In view of the computations (\ref{III.c46}), (\ref{III.c47}), (\ref{III.c48}), since $\bp\in W^{1,p}(D^2)$ for every
$p<+\infty$, we have the existence of $\bk$ in $L^{2,1}(D^2)$ such that
\be
\label{III.68}
\lf\{
\begin{array}{l}
\ds{\mathcal{L}}_{\bn}\bw=\bk\quad\quad\mbox{ in }D^2\\[5mm]
\ds\bw=0\quad\quad\mbox{ on }\p D^2\quad.
\end{array}
\rg.
\ee
We proceed now to the following Hodge decomposition of $\nabla \bw-3\pi_{\bn}(\nabla\bw)$
on $D^2$ : $\nabla \bw-3\pi_{\bn}(\nabla\bw)=\nabla \bC+\nabla^\perp \bD$ with the boundary
conditions $C=0$ on $\p D^2$ and $\p D/\p\nu=0$ on $D^2$. The following equations then holds
In one hand
\be
\label{III.69}
\lf\{
\begin{array}{l}
\ds\Delta \bC=div(\bw\wedge\nabla^\perp\bn)+\bk\quad\quad\mbox{ in }D^2\\[5mm]
\ds \bC=0\quad\quad\mbox{ on }\p D^2
\end{array}
\rg.
\ee
In the other hand
\be
\label{III.69b}
\lf\{
\begin{array}{l}
\ds\Delta \bD=3\, div(\pi_{\bn}(\nabla^\perp \bw))\quad\quad\mbox{ in }D^2\\[5mm]
\ds\frac{\p \bD}{\p \nu}=0\quad\quad\mbox{ on }\p D^2
\end{array}
\rg.
\ee
Using now the fact that the space of $L^2$ functions on $D^2$ having first derivatives in $L^{2,1}$ embedds in $L^\infty$, we get the a-priori estimates
\be
\label{III.70}
\begin{array}{l}
\ds\|\nabla \bC\|_{L^\infty(D^2)}+\|\nabla^2 \bC\|_{L^{2,1}(D^2)}\le C\|\Delta \bC\|_{L^{2,1}(D^2)}\\[5mm]
 \quad\quad \ds \le C\ \lf[\|k\|_{L^{2,1}(D^2)}+\|\nabla\bn\|_{L^{2,1}(D^2)}\ \|\nabla \bw\|_{L^\infty}\rg]
\end{array}
\ee
and
\be
\label{III.71}
\begin{array}{l}
\ds\|\nabla \bD\|_{L^\infty(D^2)}+\|\nabla^2 \bD\|_{L^{2,1}(D^2)} \le C\|\Delta \bD\|_{L^{2,1}(D^2)}\\[5mm]
 \quad\quad\ds\le C\ \|\nabla\bn\|_{L^{2,1}(D^2)}\ \|\nabla \bw\|_{L^\infty}
 \end{array}
\ee
where we have used the fact that $|div(\pi_{\bn}(\nabla\bw))|\le C |\nabla \bn|\ |\nabla\bw|$.
Thus for $\ep$ having been choosed small enough in (\ref{III.67}), following the construction of lemma~\ref{lm-III.2},
we get the existence of a lipschitz solution of (\ref{III.68}) with second derivatives in $L^{2,1}$. From the uniqueness
of lemma~\ref{lm-III.1} we obtain that this solution coincides with $\bw$ and we have then established that $\bp(\rho x+x_0)\ \chi(x)$
is Lipschitz with second derivatives in $L^{2,1}$.

A similar argument including the boundary condition $\bp=0$ could be carried out for any point $x_0$ on
 the boundary of $D^2$. We then obtain that $\bp$ is Lipschitz with second derivatives in $L^{2,1}$
and lemma~\ref{lma-A.7} is proved.\hfill$\Box$
\begin{Lma}
\label{lma-A.8}
There exists $\ep_0>0$ such that for every $0<\ep<\ep_0$ the following holds. 
Let $\bn$ be a $W^{1,2}$ map from $D^2$ into the space of unit $m-2-$vectors in ${\R}^m$ satisfying
\be
\label{III.72}
\int_{D^2}|\nabla \bn|^2\ dx\le\ep\quad.
\ee
Let $\bv$ be a function in $L^2(D^2)$ such that $\nabla \bv$ is a sum of a compactly supported distribution
in the open disk and a function in $L^{2,\infty}(D^2)$ (in such a way that the trace of $\bv$ on $\p D^2$ is well defined).
Assume $\nabla\bn$ is in the Lorentz space $L^{2,1}(D^2)$ and that $\bv$ satisfies in the distributional sense the system
\be
\label{III.73}
\lf\{
\begin{array}{l}
\ds\Delta\bv-3\ div(\pi_{\bn}(\nabla \bv))-div\lf(\bv\wedge\nabla^\perp\bn\rg)=0\quad\quad\mbox{ in } D^2
\\[5mm]
\ds\bv=0\quad\quad\mbox{on }\p D^2\quad.
\end{array}
\rg.
\ee
Then $\bv$ is identically 0 in $D^2$.\hfill$\Box$
 \end{Lma}
\noindent{\bf Proof of lemma~\ref{lma-A.8}.}
We consider a smoothing $\bv_\delta$ of $\bv$ obtained by making convolution with a function with  support shrinking as one get
close to the boundary in order to ensure that $\nabla\bv_\delta\in L^{2,\infty}(D^2)$, that $\bv_\delta=\bv$ in a neighborhood
of $\p D^2$ and that $\bv_\delta\rightarrow \bv$ in $L^2(D^2,{\R}^m)$. Let $\bp_{i}$ be a sequence of normalised eigenvectors
of ${\mathcal {L}}_{\bn}$ in $W^{1,2}_0(D^2,{\R}^m)$ and forming an Hilbert orthonormal Basis of $L^2(D^2,{\R}^m)$
(the existence of such a Basis is obtained by combining the result of Lemma~\ref{lm-III.1} and Hilbert-Schmidt theorem).
Denote $\la_i$ the corresponding eigenvalues. From lemma~\ref{lm-III.1} again we know that $\la_i\ne 0$.
We have that
\be
\label{III.74}
\int_{D^2}\bv_\delta\cdot\bp_i=\la_i^{-1}\int_{D^2}\bv_\delta\cdot{\mathcal{L}}_{\bn}\bp_i\quad.
\ee
From lemma~\ref{lma-A.7} $\bp_i$ is Lipschitz and $\nabla^2\bp_i\in L^{2,1}$. Therefore, since
$\nabla \bv_\delta\in L^{2,\infty}$ and since both $\bp_i$ and $\bv_\delta$ are 0 on $\p D^2$, we have clearly that for any $p>2$
 \be
 \label{III.75}
 \int_{D^2}\bv_\delta\cdot{\mathcal{L}}_{\bn}\bp_i=\lf<{\mathcal L}_{\bn}\bv_\delta,\bp_i\rg>_{W^{-1,p'},W^{1,p}_0}
 \ee
 Observe now that $\Delta\bv_\delta$ converges to $\Delta\bv$ in $H^{-2}(D^2)$ which is dual to $W^{2,2}_0(D^2)$,
 that $div(\pi_{\bn}(\nabla \bv_\delta))$ converges to $div(\pi_{\bn}(\nabla \bv))$ in $W^{-1,1}\oplus H^{-2}(D^2)$
 and that $div(\bv_\delta\wedge\nabla^\perp\bn)$ converges to $W^{-1,1}(D^2)$. So using that $\bp_i$ is in
  $W^{1,\infty}\cap W^{2,2}( D^2)$ and interpreting the duality in the right-hand-side of (\ref{III.75}) by the mean of these
  norms we can pass in the limit as $\delta$ goes to zero and we have that
\be
\label{III.76}
\int_{D^2}\bv\cdot{\mathcal{L}}_{\bn}\bp_i =\lf<{\mathcal L}_{\bn}( \bv),\bp_i\rg>_{H^{-2}\oplus W^{-1,1},W^{2,2}_0\cap W^{1,\infty}_0}
\ee
Combining (\ref{III.73}), (\ref{III.74}) and (\ref{III.76}) we obtain that for every $i$ $\int_{D^2}\bv\cdot\bp_i=0$ and since $(\bp_i)$
realizes an Hilbert base of $L^2$ we deduce that $\bv$ is identically zero on $L^2$ which ends the proof of lemma~\ref{lma-A.8}.
\hfill$\Box$

\begin{Lma}
\label{lma-A.9}
There exists $\ep_0>0$ such that for every $0<\ep<\ep_0$ , there exists a constant $C>0$ independent of $\ep$
 such that the following holds. Let $\bn$ be a lipschitz map from $D^2$ into the space of
 unit $m-2-$vectors in ${\R}^m$ satisfying
\be
\label{III.63a}
\int_{D^2}|\nabla \bn|^2\ dx\le\ep\quad,
\ee
Let $\bg$ be a map in $L^p(D^2)$ for some $p>1$. Let $\bv$ be an $L^2$ solution  of 
\be
\label{III.64a}
\Delta\bv-3\ div(\pi_{\bn}(\nabla \bv))-div\lf(\bv\wedge\nabla^\perp\bn\rg)=0\quad\quad\mbox{ in } D^2\quad.
\ee
Then we have
\be
\label{III.65a}
\|\nabla^2 \bv\|_{L^p(D^2_{1/2})}\le C\lf[\| \bg\|_{L^p(D^2)}+\|\bv\|_{L^2(D^2)}\rg]\quad,
\ee
where $D^2_{1/2}$ is the disk centered at $0$ with radius $1/2$.\hfill$\Box$
\end{Lma}

The proof of this result goes following the same line as the one followed in some of the previous lemma.
\medskip

\noindent{\it Remark:} The assumption that $\bv$ is in $L^2$ is not optimal in the previous lemma
for the result to be true. This is however sufficient for the use we are making of it.


\begin{thebibliography}{99}
\bibitem[Bla]{Bla} Blaschke, Wilhelm "Vorlesungen \"uber differential geometry" III, Springer (1929).
\bibitem[Bry]{Bry} Bryant, Robert L. "A duality theorem for Willmore surfaces." 
J. Differential Geom. 20 (1984), no. 1, 23--53. 
\bibitem[Che]{Che} Chen, Bang-yen "Some conformal invariants of submanifolds and their applications."  Boll. Un. Mat. Ital. (4) 10 (1974), 380--385.
\bibitem[Fe]{Fe} Federer, Herbert  "Geometric measure theory", Springer (1969).
\bibitem[FJM]{FJM} Friesecke, Gero; James, Richard D.; MŸller, Stefan "A theorem on geometric rigidity and the derivation of nonlinear plate theory from three-dimensional elasticity". Comm. Pure Appl. Math. 55 (2002), no. 11, 1461--1506.
\bibitem[Haw]{Haw} S. W. Hawking, Gravitational radiation in an expanding universe, J. Math. Phys. 9 (1968), 598--604.
\bibitem[Hel]{Hel} H\'elein, Fr\'ed\'eric "Harmonic maps, conservation laws and moving frames". 
 Cambridge Tracts in Mathematics, 150. Cambridge University Press, Cambridge, 2002.
 \bibitem[Hef]{Hef} Helfrich, W.  Z. Naturforsch., C28 (1973), 693-703.
 \bibitem[Hub]{Hub} Huber, Alfred  "On subharmonic functions and differential geometry in the large."
  Comment. Math. Helv. 32 1957 13--72.
\bibitem[HI]{HI} Huisken, G.; Ilmanen, T. "The Riemannian Penrose inequality." Internat. Math. Res. Notices 1997, no. 20, 1045--1058.
\bibitem[KS1]{KS1}Kuwert, Ernst; Sch\"atzle, Reiner "The Willmore flow with small initial energy". J. Differential Geom. 57 (2001), no. 3, 409--441.
\bibitem[KS2]{KS2} Kuwert, Ernst; Sch\"atzle, Reiner "Gradient flow for the Willmore functional". Comm. Anal. Geom. 10 (2002), no. 2, 307--339.
\bibitem[KS3]{KS3} Kuwert, Ernst; Sch\"atzle, Reiner "Removability of point singularities of Willmore surfaces". Ann. of Math. (2) 160
  (2004), no. 1, 315--357.
\bibitem[LY]{LY}  Li, Peter; Yau, Shing Tung "A new conformal invariant and its applications to the Willmore conjecture and the first eigenvalue of compact surfaces." Invent. Math. 69 (1982), no. 2, 269--291.
\bibitem[Ri1]{Ri1} Rivi\`ere, Tristan "Conservation laws for conformally invariant variational problems" to appear in Invent. Math. (2006).
\bibitem[Mon]{Mon} Montiel, Sebasti\`an "Willmore two-spheres in the four-sphere." 
Trans. Amer. Math. Soc. 352 (2000), no. 10, 4469--4486. 
\bibitem[MS]{MS}  MŸller, S.;  Sver‡k, V. "On surfaces of finite total curvature." J. Differential Geom. 42 (1995), no. 2, 229--258.
\bibitem[Si1]{Si1} Simon, Leon "Lectures on geometric measure theory." Proceedings of the Centre for Mathematical Analysis, Australian National University, 3. Australian National University, Centre for Mathematical Analysis, Canberra, 1983.
\bibitem[Si2]{Si2} Simon, Leon "Existence of surfaces minimizing the Willmore functional." Comm. Anal. Geom. 1 (1993), no. 2, 281--326.
\bibitem[Sim]{Sim} Simonett, Gieri "The Willmore flow near spheres." Differential Integral Equations 14 (2001), no. 8, 1005--1014. 
\bibitem[Ta1]{Ta1} Tartar, Luc "Remarks on oscillations and Stokes' equation.  Macroscopic modelling of turbulent flows" (Nice, 1984), 24--31, Lecture Notes in Phys., 230, Springer, Berlin, 1985.
\bibitem[Ta2]{Ta2} Tartar, Luc "Imbedding theorems of Sobolev spaces into Lorentz spaces." 
Boll. Unione Mat. Ital. Sez. B Artic. Ric. Mat. (8) 1 (1998), no. 3, 479--500. 
\bibitem[Tho]{Tho} Thomsen, G. " \"Uber konforme Geometrie, I" Grundlagen der Konformen Fl\"achentheorie. Abh Math. Sem. Hamburg, (1923) 31-56. 
\bibitem[Wei]{Wei} Weiner, Joel ``On a problem of Chen, Willmore, et al'' Indiana Univ Math. J, 27, no 1 (1978), 19-35.
\bibitem[Wen]{Wen} Wente, Henry C. "An existence theorem for surfaces of constant mean curvature." J. Math. Anal. Appl. 26 1969 318--344.
\bibitem[Whi]{Whi} White, James H.  "A global invariant of conformal mappings in space". 
Proc. Amer. Math. Soc. 38 (1973), 162--164. 
\bibitem[Wil]{Wil} Willmore, T. J. "Note on embedded surfaces." Ann. Stiint. Univ. "Al. I. Cuza" Iasi. Sect. I a Mat. (N.S.) 11B 1965 493--496.
 \end{thebibliography}
\end{document}